\documentclass[11pt]{article}

\usepackage{amsfonts,amsthm,amssymb}
\usepackage[cp866]{inputenc}

\usepackage[T2A]{fontenc}

\usepackage{epsfig}
\usepackage{graphicx}

\begin{document}

\newtheorem{lm}{Lemma}
\newtheorem{theorem}{Theorem}
\newtheorem{prop}{Proposition}
\newtheorem{df}{Definition}
\newtheorem{remark}{Remark}
\newtheorem{corollary}{Corollary}
\newtheorem{ex}{Example}

%\iffalse

%\begin{document}
\begin{center}{\large\bf Elements of contemporary mathematical theory of dynamical chaos.\\ Part 1. Pseudohyperbolic attractors.}
\end{center}
\medskip
\begin{center}{\bf Gonchenko S.V.$^1$, Gonchenko A.S.$^1$,  Kazakov A.O.$^{2,1}$, Kozlov A.D.$^1$ }
\end{center}
\medskip
\begin{center}{$^1$ Lobachevsky state university of Nizhny Novgorod, Russia}
\end{center}
\begin{center}{$^2$ High School of Economics, Nizhny Novgorod, Russia}
\end{center}

\begin{center} {e-mail: sergey.gonchenko@mail.ru, phone:+79065567766; e-mail: agonchenko@mail.ru, phone:+79527617397;  e-mail: kazakovdz@yandex.ru, phone: +79081659694; e-mail: kozzzloff@list.ru, phone: +79506128267}
\end{center}
\bigskip

\begin{abstract}
The paper deals with topical issues of modern mathematical theory of dynamical chaos and its applications.
At present, it is customary to assume that dynamical chaos  in finite-dimensional smooth systems can exist in three different forms.
This is {\em dissipative chaos}, the mathematical image of which is a strange attractor; {\em conservative chaos}, for which the entire phase space is a large ``chaotic sea`` with randomly spaced elliptical islands inside it; and {\em mixed dynamics}, characterized by the principal inseparability in the phase space of attractors, repellers and conservative elements of dynamics. In the present paper (which opens a cycle of three our papers), elements of the theory of pseudo-hyperbolic attractors of multidimensional maps are presented. Such attractors, as well as hyperbolic ones, are genuine strange attractors, but they allow the existence of homoclinic tangencies. We give a mathematical definition of a pseudo-hyperbolic attractor for the case of multidimensional maps, from which we derive the necessary conditions for its existence in the three-dimensional case, formulated using the Lyapunov exponents. We also describe some phenomenological scenarios for the appearance of pseudo-hyperbolic attractors of various types in one-parameter families of three-dimensional diffeomorphisms, we propose new methods for studying such attractors (in particular, a method of saddle charts and a modified  method of Lyapunov diagrams). We consider also
%orientable and non-orientable
three-dimensional generalized H\'enon maps as examples.

In the second part, we plan to review the theory of spiral attractors, which compose  an important class of attractors often meeting in applications. The third part will be devoted to mixed dynamics -- a new type of chaos, which is characteristic, in particular for reversible systems, i.e. systems invariant with respect to the time reversal. It is well known that such systems are met in many problems of mechanics, electrodynamics and other fields of natural science.

\end{abstract}

\section*{Introduction}

At present, one can distinguish three relatively independent and different forms of dynamic chaos of smooth finite-dimensional systems -- ``dissipative chaos'', ``conservative chaos'' and ``mixed dynamics''. The dissipative chaos is characterized by the existence of {\em a strange attractor} -- a nontrivial attracting closed invariant set which is located in the phase space of a system inside some absorbing domain, into which all orbits intersecting its boundary enter. Unlike the dissipative chaos, conservative chaos spreads over the whole phase space.
Mixed dynamics is a new third type of dynamical chaos which is characterized by the fact that stable elements of dynamics (attractors) coexist with completely unstable ones (repellers), and, moreover, they are principally inseparable from each other and from conservative elements of dynamics (so-called reversible cores), see more details in \cite{G16,GT17}.

Recall that the Conley's theorem \cite{Co78}
implies that
any smooth system on a compact manifold $M$ has an attractor ${\cal A}$ and a repeller ${\cal R}$ (i.e. attractor at time reversal). Then we have always that ${\cal A}\cap {\cal R} = \emptyset$ in the case of dissipative chaos and ${\cal A} = {\cal R} = M$ in the case of conservative chaos.
For mixed dynamics we have the third logically possible case when ${\cal A}\cap {\cal R} \neq \emptyset$ and ${\cal A} \neq {\cal R}$. This automatically implies that the existence of fourth type of chaos  is principally impossible.

The phenomenon of mixed dynamics was discovered, in fact, in the work
\cite{GST97}, where, in particular, it was proved that, in the space of two-dimensional diffeomorphisms, there are open regions (the so-called Newhouse regions\footnote{Recall, that Newhouse regions are open regions from the space of $C^r$-smooth systems, $r\geq 2$, where systems with homoclinic tangencies are dense, i.e. systems having saddle periodic orbits whose invariant manifolds intersect  nontransversely. In turn, Newhouse regions exist in any neighbourhood of any system with homoclinic tangency \cite{N79,GST93b,PV94,Rom95}.}) in which diffeomorphisms with infinitely many stable, completely unstable and saddle periodic orbits are dense and, moreover, in this case the closures of sets of orbits of different types have nonempty intersections. Naturally, for any definition of  attractor, this must be, in any case, a stable closed invariant set that should %
contain all stable periodic orbits, if they exist. The same fact (for completely unstable periodic orbits) must also hold
%be satisfied
for repellers. Thus, in \cite{GST97} it was shown that attractor can intersect with repeller.
The mathematical justification of this phenomenon (of attractor and repeller merger)  is given quite recently, see e.g. \cite{G16,GT17}.
%, which formally agrees with the Conley  theorem.

As for the strange attractors to which this article is devoted, their generally accepted definition, which would be suitable for all occasions,  does not exist up to now. The only exceptions are the so-called {\em genuine strange attractors} whose definition includes two main conditions: 1) the existence of an absorbing region in the phase space, where the attractor exists,
 %(in which enter all orbits crossing the boundary of this region),
and 2) the instability of orbits of the attractor, which means that each orbit of the attractor has a positive maximal Lyapunov exponent. It is also assumed that properties 1) and 2) are satisfied for all nearby systems. On the other hand, the so-called {\em quasiattractors} are also referred  to the strange attractors, with due reason, see discussions in \cite{AfrSh83a,Sh97,GST97a}. Here, under the quasiattractors we mean a non-trivial attracting invariant set that either contains stable periodic orbits of very large periods (and with very narrow domains of attraction), or such orbits appear under arbitrarily small smooth perturbations. This is connected with the fact that the quasiattractors admit the existence of saddle periodic orbits whose stable and unstable invariant manifolds intersect nontransversely. In turn, bifurcations of such homoclinic tangencies, when certain conditions-criterion are fulfilled
\cite{GaS73,G83,GST93c,GST96,GST08}, lead to the birth of asymptotically stable periodic orbits, attracting invariant tori, small H\'enon-like attractors, \cite{MV93,PV94,GStT96,Colli,Hom02}, and even small Lorenz-like attractors, \cite{GMO06,GST09,GO13,GOT14,GO17}, and so on.

We also note that, in dissipative systems with special structures, strange attractors of other types may exist, which do not formally fit into this scheme.
For example, nonsmooth or discontinuous systems may have attractors with a singularly hyperbolic behavior of orbits, in the sense that  Lyapunov exponents are not defined for some orbits (by virtue of the nonsmoothness of the system itself), although
there are both absorbing regions and the necessary instability of orbits on the attractor. Examples of such attractors are the Lozi attractor \cite{Lozi} and the Belykh attractor \cite{Belykh}. A completely different type of complex non-periodic behavior of orbits is demonstrated by the so-called strange non-chaotic attractors that arise in special models that have the structure of a direct product of a nonspecific dynamical system and a quasi-periodic system. They are characterized by the fact that one of the Lyapunov exponents is zero for any orbit, and the remaining ones are less than zero (there is also allowed the existence of a set of zero measure of orbits with a positive Lyapunov exponent).
%there is also the existence of a small number (of zero measure) of orbits with a positive exponent).
For more details,  see e.g. \cite{Sncha}.

One can say that the most of known strange attractors of smooth dynamical systems, in particular, many such attractors in systems from applications, are, in fact, quasiattractors.
%We can say that most of the known strange attractors of smooth dynamical systems, including many of them from %those encountered in
%applications, are, in essence, quasi-attractors.
Examples of such attractors are the numerous `torus-chaos' attractors arising from break-down of two-dimensional tori, \cite{AfrSh83b}; attractors in Chua circuits \cite{Chua}; the H\'enon attractor  \cite{Henon76,BenCar}; attractors in periodically perturbed two-dimensional systems with a homoclinic figure-eight of a saddle \cite{GSV13} and many others.

We should also note that a special class of quasiattractors consists of the so-called {\em spiral attractors} that are related to the existence homoclinic orbits to a saddle-focus equilibrium. Such attractors are often found in applications, and examples of such attractors are well known. In particular, spiral attractors of three-dimensional flows, such as the Roessler attractor \cite{Rossler1976,Rossler1976_2}, attractors in the Arneodo-Coullet-Tresser  models \cite{ArnCoulTres1980,ArnCoulTres1981,ArnCoulTres1982} (called also as ACT-attractors), in the Rosensweig-MacArthur models \cite{RosMac} etc.
It is interesting that all these attractors appear in flows according to a fairly simple and universal phenomenological scenario proposed by Shilnikov \cite{Sh86}.
The main specific of spiral attractors of such type is that they  contain, at certain values of parameters, a saddle-focus with the {\em two-dimensional}  unstable invariant manifold.
With some natural modifications, the Shilnikov scenario was also transferred to the case of three-dimensional maps \cite{GGS12}.
Therefore, for such spiral attractors, we proposed in \cite{GGKT14} the generalizing title:  ``Shilnikov attractor`` for flows (these include, in particular, the  mentioned above  Roesler attractor and ACT-attractors); or  ``Shilnikov discrete attractor'' for maps. Various examples of such discrete attractors were found in  three-dimensional H\'enon maps \cite{GGS12,GGKT14,GGOT13,GG16,GK16} for which very interesting results were obtained, and we plan to review them in the next part of our paper.

Until recently, only the hyperbolic attractors and Lorenz attractors could be considered as  genuine strange attractors of smooth dynamical systems. However, the situation has been changed after the work by Turaev and Shilnikov \cite{TS98}, where a new class of genuine strange attractors was introduced, the so-called {\em wild hyperbolic attractors}. These attractors, unlike hyperbolic and Lorenz ones, admit the existence of homoclinic tangencies, but they do not contain stable periodic orbits and any other stable invariant subsets that do not arise also for small smooth perturbations. Systems with wild hyperbolic attractors belong to Newhouse domains.
%, i.e. open (in $C^2$-topology) domains from the space of dynamical systems in which systems with homoclinic %tangencies are dense.
\footnote {The term ``wild'' goes back to the Newhouse paper \cite{N79}, in which the concept of ``wild hyperbolic set'' was introduced, i.e. such a uniformly hyperbolic invariant set, in which, among its stable and unstable invariant manifolds, there are  always those that intersect nontransversely, and this property is preserved for all small $C^2$-smooth perturbations.} However, these tangencies, unlike homoclinic tangencies in systems with quasiattractors, do not lead to the appearance of stable periodic orbits \cite{GST93c,GST96,GST08}, see also Sec.~\ref{sec:phomtan}.

In \cite{TS98}, an example of a four-dimensional flow with a wild spiral attractor containing an equilibrium state of the saddle-focus type was also constructed. One of the main features of the Turaev-Shilnikov spiral attractor is that it possesses an {\em pseudo-hyperbolic structure}. Speaking shortly,  this feature means that, in a neighborhood $D$ of the attractor (one can think that $D$ is some its absorbing region), there is a ``weak'' version of hyperbolicity: there is a partition of $D$ into  subsets (strongly contracting and volume expanding), that are transversal and invariant with respect to the differential, such that on one of them there takes place exponential contraction along all directions, and on the other -- exponential expansion of the volume. It is also required that such a partition depends continuously on a point from $ D $; the corresponding coefficients of compression and extension, and also the angles between the tangent vectors of the subspaces are uniformly bounded; in the subspace, where the volume is expanded,  if there are any contractions, then all of them are uniformly weaker than any contraction in a strongly contracting space.

We note that the pseudohyperbolicity conditions are verified for points of the absorbing domain $D$. If they are fulfilled, then, as shown in \cite{TS98}, the attractor exists and it is unique,  each of its trajectories has a positive maximal Lyapunov exponent (this follows from the property of expanding volumes)\footnote{In this case  the attractor is, in fact,  the Ruelle attractor  \cite {Ruelle}, i.e. closed, invariant, (asymptotically) stable and chain-transitive set, for more details see \cite{GT17} and Sections~\ref{sec:fen}.}.

In fact, in the paper \cite{TS98} it was laid the foundations of a very promising theory of pseudohyperbolic strange attractors. New examples of such attractors were also found shortly. Thus, in the paper \cite{GOST05} it was shown that for three-dimensional H\'enon maps of the form
\begin{equation}
\bar x = y,\;\bar y = z,\; \bar z = M_1 + B x + M_2 y - z^2,
\label{3DHM1}
\end{equation}
where $M_1,M_2,B$ are parameters ($B$ is the Jacobian), in a certain domain of  parameter values adjoining the point $A^* = (M_1 = 1/4, M_2 = 1, B = 1) $, there exist discrete Lorenz attractors.

The pseudohyperbolicity of such attractors was claimed in \cite{GOST05}  on the basis of the fact that for values of parameters close to $A^*$, the second power of map (\ref{3DHM1}) in some neighborhood of the saddle fixed point can be represented as the Poincar\'e map of a periodically perturbed Shimizu-Morioka system, which, in turn, has the  Lorenz attractor  \cite{ASh86,ASh93}. If the perturbation is sufficiently small (which is determined by the closeness of the values of parameters to $A^*$), then the desired pseudohyperbolicity should naturally be inherited from the pseudohyperbolicity of the Lorenz attractor \cite{Tucker99,TurOvs17}. In particular, in \cite{TS08} it was shown that the property of pseudohyperbolicity of flows is also preserved for their Poincar\'e maps for small periodic perturbations.

\begin{figure} [tb]
  \centerline{
  \includegraphics[width=16cm]{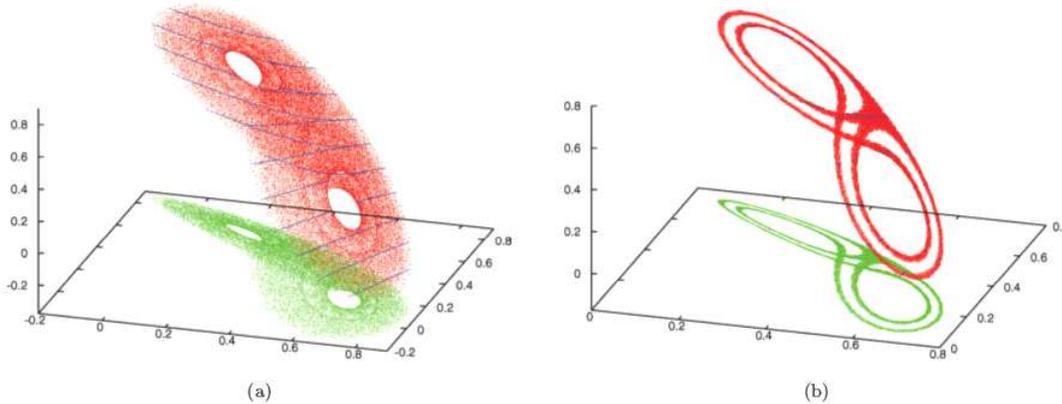}}
  \caption{{\footnotesize Portraits of discrete Lorenz attractors in the case of  map (\ref{3DHM1}) with (a) $ M_1 = 0, M_2 = 0.85, B = 0.7 $ and (b) $ M_1 = 0, M_2 = 0.825, B = 0.7 $. In both cases, about  $10^5$ iterations of a single initial point on the attractor is shown. The projections of attractors onto the plane $(x, y)$ and some slices of the attractor by the plane $z = \mbox{const}$ are also shown (although these slices look like lines, in fact they have a complex Cantor structure).}}
     \label{ExDLA1}
\end{figure}

In Fig.~\ref{ExDLA1}  we show examples of discrete Lorenz attractors for map (\ref{3DHM1}). We note that the phase portraits of these attractors are very similar to portraits of the Lorenz attractors for flow. However, we see that the corresponding values of the parameters ($ M_1 = 0, M_2 = 0.85, B = 0.7 $ in the case of Fig.~\ref{ExDLA1}(a) and $ M_1 = 0, M_2 = 0.825, B = 0.7 $ in the case of Fig.~\ref{ExDLA1}(b)) are not nearly close to $A^*$. Therefore, the conditions for pseudohyperbolicity of such attractors need to be checked additionally.

Such a task seems very complicated. In fact, we can only test some necessary conditions. For example, for a pseudohyperbolic attractor of a three-dimensional map, its Lyapunov exponents $\Lambda_1> \Lambda_2> \Lambda_3$ must satisfy such conditions
\begin{equation}
\Lambda_1 >0,\; \Lambda_1 + \Lambda_2 > 0,\;  \Lambda_1 + \Lambda_2 + \Lambda_3<0.
\label{lyapcond}
\end{equation}
The first and third conditions indicate that the observed attractor is strange, and the second one means that there is a expansion of two-dimensional areas on it.

However, the Lyapunov exponents are certain average characteristics of the orbits on an attractor, therefore, in principle, it is possible that the attractor has very small ``windows'' (whose sizes may be less than any reasonable accuracy of calculations), where conditions (\ref{lyapcond}) for the corresponding orbits are violated.
At our request, conditions for  pseudohyperbolicity of the attractors from Fig.~\ref{ExDLA1} were verified (by means of the interval arithmetic methods \cite{Tucker99}), by mathematicians from the University of Uppsala, Sweden, J. Figueros and W. Tucker, who have obtained very interesting and fine results.
Namely, they have shown that a stable periodic orbit with extremely small basin of attraction exists inside the attractor from Fig.~\ref{ExDLA1}a, while,  the attractor from Fig.~\ref{ExDLA1}b is a genuine pseudohyperbolic attractor. 
Besides, independently and by other methods, similar results were obtained by
%Saratov mathematicians
S.P. Kuznetsov and P.V. Kuptsov, who have also verified pseudohyperbolicity of some other attractors from our paper  \cite{GG16}. We hope that these interesting results will be published in the near future.

The existence of homoclinic tangencies naturally kills the hyperbolicity, however, in general, the pseudo-hyperbolicity is not violated. For example, both attractors from Fig.~\ref{ExDLA1} contain a saddle fixed point with multipliers
$\lambda_1,\lambda_2,\lambda_3$ such that $\lambda_1< -1,0<\lambda_2<1,-1<\lambda_3<0$ and  $|\lambda_2|>|\lambda_3|$, where $\lambda_1\lambda_2\lambda_3 = B = 0.7<1$ and the saddle value $\sigma = |\lambda_1\lambda_2|$ is greater than 1. Then the homoclinic tangencies, that inevitably occur here,
will be typically such as in Fig. 3, and bifurcations of such tangencies\footnote{In the general case, such (quadratic) homoclinic tangencies are called simple \cite{GST08}, and in the case $\sigma> 1 $ they do not destroy the pseudohyperbolicity (if the fixed point itself is pseudohyperbolic), although they are certain indicators of  wild hyperbolicity (for more detail see \cite{GST93b,GST93c,GST08}).
}
do not lead to the birth of stable periodic orbits \cite{GST93c,GST96,GST08}.

\begin{remark} {\rm
On the other hand, stable periodic orbits are necessarily born if $\sigma <1 $, or if a fixed (periodic) point of an attractor is a saddle-focus (whether with one-dimensional or with two-dimensional unstable manifold). In particular, the spiral attractors of three-dimensional smooth maps or flows are always quasi-attractors. In this connection, the following problem seems to be very interesting: {\em let a three-dimensional diffeomorphism have a strange attractor containing a saddle fixed point with the two-dimensional unstable invariant manifold, then this attractor is the quasi-attractor.}\footnote{This problem seems to be very difficult, and its solution is connected, for example, with the proof of the existence of the so-called non-simple homoclinic tangencies \cite{Tat01,GGTat07,GOT14}, examples of which are shown in Fig.~\ref{htnonsimple} -- only here the direction of the arrows should be reversed, so that the unstable manifold of the point $O$ becomes two-dimensional. We note that bifurcations of such tangencies lead to the birth of stable periodic orbits \cite{Tat01}. In turn, the appearance of non-simple tangencies in the case under consideration is to be very expected, because the two-dimensional unstable manifold should fold  infinitely many times  in different directions (it is as if we tried to ``package'' the two-dimensional plane into a three-dimensional cube, while avoiding the appearance of sharp corners).}
}
\label{rem01}
\end{remark}

For this reason, in the present paper we consider only such strange attractors of three-dimensional maps that contain fixed points of saddle type with one-dimensional unstable manifolds and with saddle value $\sigma$ greater than 1. Moreover, the main attention we pay to the so-called  homoclinic attractors that contain {\em exactly one} saddle fixed point. As shown in the works \cite{GOST05,GGS12,GGKT14,GG16}, this direction is rather  promising.

The content of the paper. In Sec.~\ref{sec:phomtan}  we discuss  basic properties of pseudohyperbolic maps (diffeomorphisms) and types of homoclinic tangencies that support or break pseudohyperbolicity. In Sec.~\ref{sec:fen}, phenomenological scenarios of the appearance of strange homoclinic attractors in one parameter families of three-dimensional diffeomorphisms
%, both orientable and nonorientable,
are discussed. In Sec.~\ref{Ex8Lor} we give examples of such attractors in the case of three-dimensional generalized H\'enon maps. The definition of a pseudohyperbolic diffeomorphism is given in the Appendix, Sec.~\ref{App1}.

\section{Pseudohyperbolicity and homoclinic tangencies} \label{sec:phomtan}

In this section we consider the basic concepts of the theory of pseudohyperbolic strange attractors. The notion of pseudohyperbolic system  was introduced by Turaev and Shilnikov: the flow case was considered in   the work \cite{TS98} and this notion for discrete dynamical systems, diffeomorphisms,  was given in \cite{TS08}, see also  \cite{Sat10,GGK13} and Appendix (Sec.~\ref{App1}) to this article.
In short, the pseudo-hyperbolicity of a diffeomorphism
$f$ on some domain ${\cal D}$ means the following:
\begin{itemize}
\item[(i)]
at each point of ${\cal D}$, there are two transversal linear subspaces $N_1$ and $N_2$ that are continuously dependent on the point and invariant with respect to the differential
$Df$ of $f$;
\item[(ii)]
$Df$ is exponentially strongly contracting on $N_1$ and expanding (exponentially) volumes on $N_2$ (here the word `` strongly '' means that any possible contraction in $N_2$ is uniformly weaker than any contraction in $N_1$);
\item[(iii)]
angles between $N_1$ and $N_2$ are uniformly separated from 0 as well as exponents of contraction in $N_1$ and of expansion of area in $N_2$ are uniformly bounded and separated from 0 (see Def.~\ref{def:psevhyp} from Sec.~\ref{App1}).
\end{itemize}

Thus, unlike hyperbolicity, here is not required the existence of uniform expansions in
$N_2$ along all directions. Nevertheless, the pseudohyperbolicity, the same as  the hyperbolicity, is preserved under small smooth perturbations \cite{TS98,TS08}. Therefore, if the (pseudohyperbolic) diffeomorphism $f$ has an attractor inside ${\cal D}$, then this attractor is strange, since the expansion of volumes in $N_2$ guarantees the existence of a positive maximal Lyapunov exponent for any orbit. In other words, pseudohyperbolic attractors are genuine attractors.

However, in contrast to the hyperbolic and Lorenz attractors,
%{\em homoclinic tangencies} can exist in
pseudo-hyperbolic attractors can possess {\em homoclinic tangencies}.
Moreover, if it is not known in advance that the strange attractor is hyperbolic, then in addition to transverse homoclinic orbits (at points of which stable and unstable invariant manifolds of saddle periodic orbits intersect transversally), there should also exist nontransversal ones -- homoclinic tangencies.
By itself, the appearance of a homoclinic tangency is not something exceptional: this is a codimension one bifurcation phenomenon when a quadratic homoclinic tangency appear.
However, as S. Newhouse has shown \cite{N79}, this single bifurcation can imply very complicated structure of the bifurcation set. In particular, he has proved that, arbitrarily close to any two-dimensional diffeomorphism with a homoclinic tangency, there are open regions in which diffeomorphisms with homoclinic tangencies are dense. In multidimensional case the existence of Newhouse regions near any system with a homoclinic tangency was established later, see e.g.  \cite{GST93b,PV94,Rom95}.
%These regions are called now the Newhouse regions. where   yin any neigh for a one parameter family generally %unfolding a quadratic homoclinic tangency, there are open intervals (Newhouse intervals) in which  parameter %values corresponding to the existence of (quadratic) homoclinic tangencies are dense.
Dynamics of systems from Newhouse regions is extremely rich. So, as it was established in \cite{GST93a,GST99}, in these regions there are dense systems with  infinitely many homoclinic tangencies of any orders, systems with arbitrarily degenerate periodic orbits etc.
%can occur c, these tangencies can be degenerate \cite{GST93a,GST99}, which in turn means the possibility of the %appearance of arbitrarily degenerate periodic trajectories, etc.
All this means that bifurcations of homoclinic tangencies can not be studied completely, for example, by means of finite-parameter families -- the traditional apparatus of the classical theory of bifurcations, see more discussions in \cite{GST91}. Here, of necessity, problems of a different kind arise, connected with the investigation of the basic bifurcations and the basic characteristic properties of such systems. Moreover, what is very important and interesting, the question of which bifurcations and which characteristic properties are the main ones should be decided by the researcher himself.

In the theory of strange attractors of smooth dynamical systems, one of the most important problems relates to  an identification whether a given attractor is the quasiattractor or the genuine attractor (in particular, pseudo-hyperbolic one). Sometimes we can easily identify that the attractor under consideration is a quasiattractor. So,
in the case of strange attractors of two-dimensional diffeomorphisms (if they are not hyperbolic), bifurcations of inevitable in them homoclinic tangencies  lead to the appearance of stable periodic orbits of quite large periods \cite{GaS73}, and, accordingly, any such attractor should be considered as  quasiattractor.\footnote{This is true, for example, for the H\'enon attractors, for which stable periodic orbits arise under arbitrarily small perturbations, although they may be absent (for parameter values forming a nowhere dense set of positive measure, according to the Benedics-Carleson theory \cite{BenCar}).}

\begin{figure} [tb]
  \centerline{
  \includegraphics[width=14cm]{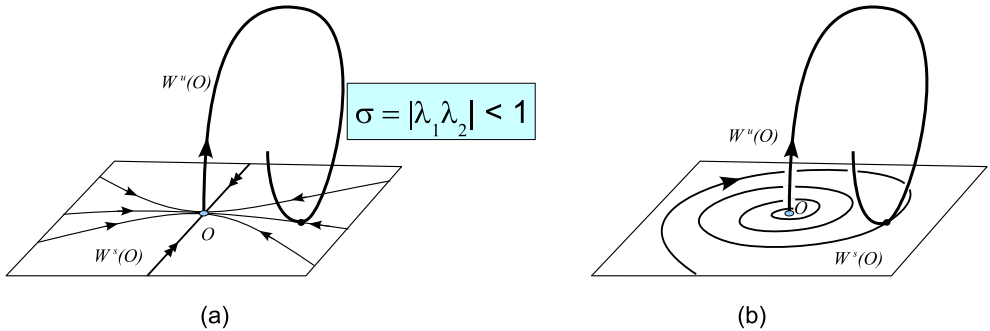}}
  \caption{{\footnotesize Examples of homoclinic tangencies whose bifurcations lead to the birth of stable periodic orbits.
  }}
   \label{homtan2}
\end{figure}

In the case of strange attractors of three-dimensional diffeomorphisms, which are one of the main subject of the present paper, the problem of an identification of their types (quasiattractor or genuine attractor) is much more complicated. However, even here, homoclinic tangencies found in attractors can be consider as peculiar indicators.
Thus, if an attractor allows homoclinic tangencies to a fixed or periodic point such as in Fig.~\ref{homtan2}, then it is definitely the quasiattractor. In the first case,
Fig.~\ref{homtan2}(a), the fixed point is a saddle with the saddle value $\sigma$ less than 1, and in the second case, Fig.~\ref{homtan2}(b), it is a saddle-focus.
The birth of stable periodic orbits under bifurcations of such homoclinic tangencies was established e.g. in \cite{GST93c,GST96,GST08} --
here it is only required that the Jacobian $J$ of the fixed point is less than one, and in the case of a saddle its unstable manifold is one-dimensional (in the case of a saddle-focus, there is no meaning whether the manifold is one-dimensional or two-dimensional).

\begin{figure} [tb]
  \centerline{
  \includegraphics[width=16cm]{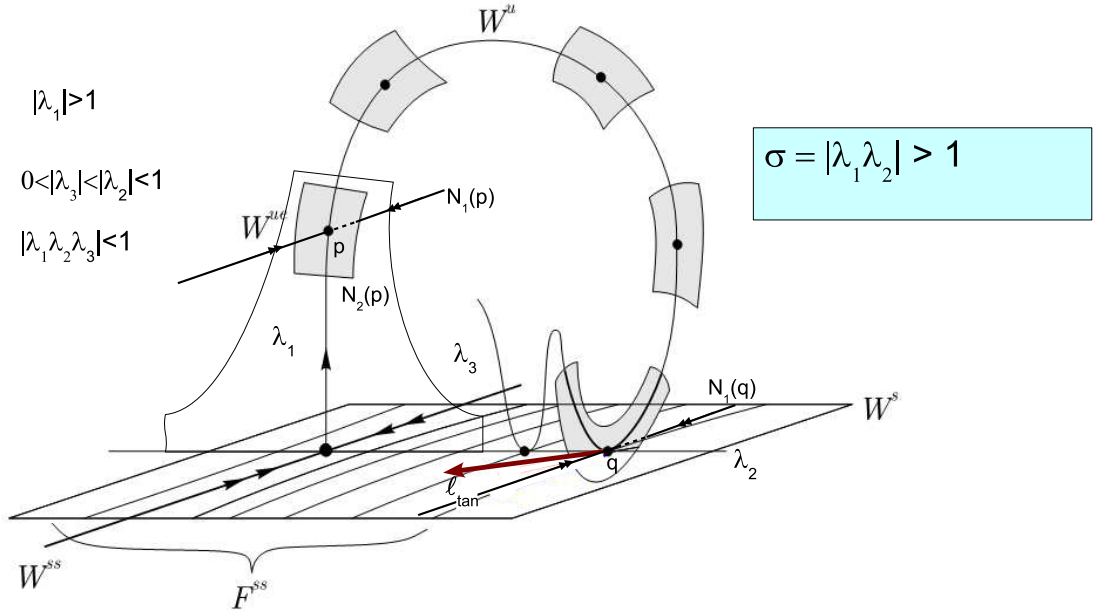}}
  \caption{{\footnotesize To the definition of simple homoclinic tangency.
  }}
   \label{htansimp1}
\end{figure}

On the other hand, it is very important that there are homoclinic tangencies that do not destroy pseudohyperbolicity.
In the case of three-dimensional diffeomorphisms, these are {\em simple homoclinic tangencies}, \cite{GST93b,GST93c}, provided that $\sigma> 1$.

Let, for example, a diffeomorphism $f$ have a saddle fixed point $O$ with real eigenvalues
$\lambda_1, \lambda_2, \lambda_3$ such that $| \lambda_1 |> 1> | \lambda_2 |> | \lambda_3 |> 0 $ under the condition that $ \sigma = | \lambda_1 | | \lambda_2 |> 1 $.
For such tangencies, the point $O$ itself is pseudo-hyperbolic: it has $N_1(O)$ as the line passing through $O$ in the direction of the eigenvector of the linearization matrix $A$ corresponding to its strong stable eigenvalue  $\lambda_3$, and $N_2(O)$ is the plane containing the eigenvectors of the matrix $A$ corresponding to the multipliers $\lambda_1$ and $\lambda_2$.
Obviously, for any point $p$ from a small neighborhood $U(O)$ of the saddle $O$ there will be the same invariant decompositions into the spaces $N_1(p)$ and $N_2(p)$. Similar invariant expansions near the entire homoclinic orbits can also be obtained if the homoclinic tangency is simple \cite{GST08}.

The simplicity of a homoclinic tangency in the case of the diffeomorphism $f$ can be defined as follows \cite{GST93c,GST08}.
Choose two any homoclinic points $p$ and $q$ in $U(O)$ such that $p\in W^u_{loc}(O)$,
$ q \in W^s_{loc}(O) $, and $f^s(p)=q$ for some integer $s$.  Define the so-called global map $T_1$ that is constructed along orbits of $f$  and acts from a small neighborhood $V(p)$ of the point $ p $ to a small neighborhood of the point $q$. Let $q = f^s(p)$ for some natural $s$, then we can write
% such that $T_1 (p) = q $ (note that $ f ^ s (p) = q $ for some natural $ s $,
$T_1 = f^s|_{V (p)}$.\footnote{Note that the local invariant manifolds $ W^u_{loc} (O) $ and $ W^s_{loc} (O) $ can always be straightened by introducing in $U(O)$ such
$C^r$-coordinates $ (x, y, z) $ in which $ W^u_{loc}(O) = \{x = 0, y = 0 \}$ and
$W^s_{loc}(O) = \{z = 0 \}$, \cite{GST08}.}
%Define the so-called. a global map $ T_1 $ that is constructed along the trajectories of the diffeomorphism under consideration and acts from a small neighborhood $ V (p) $ of the point $ p $ to a small neighborhood of the point $ q $ such that $ T_1 (p) = q $ (note that $ f ^ s (p) = q $ for some natural $ s $, then $ T_1 = f ^ s | _ {V (p)} $).
Then it is required that
\begin{itemize}
\item[$\bullet$]
the plane $DT_1(N_2 (p))$ intersects transversally with
$N_1(q)$. % and with $ W ^ s_ {loc} (O) $.
\end{itemize}
%the plane $ DT_1 (N_2 (p)) $ intersects transversally with
%$ N_1 (q) $ and with $ W ^ s_ {loc} (O) $.
Note that the curve $T_1 (W^u_{loc}(O))$ touches the two-dimensional plane
$W^s_{loc}(O)$ along the vector $\ell_{tan}$, which, in turn, has a nonzero angle with the line $ N_1 (q) $, see Fig.~\ref{htansimp1}.

\begin{figure} [tb]
  \centerline{
  \includegraphics[width=16cm]{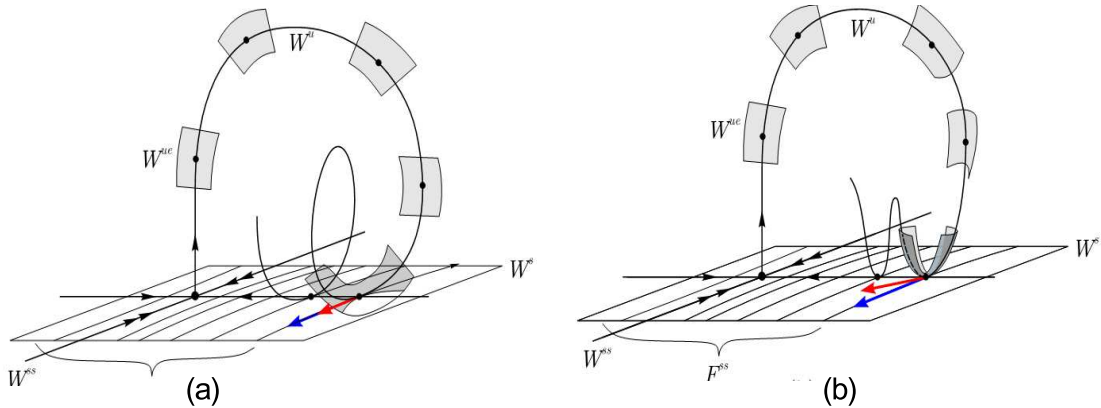}}
  \caption{{\footnotesize Two types of non-simple homoclinic tangencies: (a) when the surface  $T_1(N_2(p))$ intersects transversally with  $W^s_{loc}(O)$ but the vector $\ell_{tan}$ belongs to $N_1(q)$; (b) when the surface $T_1(N_2(p))$ touches $W^s_{loc}(O)$. }}
   \label{htnonsimple}
\end{figure}

If an attractor of a three-dimensional smooth map is pseudohyperbolic, then it can contain only simple homoclinic tangencies.\footnote{Moreover, except for quadratic tangencies, here can exist  homoclinic tangencies of arbitrarily large orders \cite{GST93a,GST99}, but they all should be simple (in the sense that,
at any homoclinic point $p$, the subspaces $N_2(p)$ and $N_1(p)$ intersect transversally, see more detail in \cite{GMLi17}.}
For any small smooth perturbations, pseudohyperbolicity is preserved \cite{TS98,TS08}. However, if these perturbations are not too small, it can  be broken.  In this case, the destruction itself can be caused by the appearance of such homoclinic tangencies as in Fig.~\ref{homtan2} (for example, the fixed point, initially with $\sigma> 1 $, in the process of evolution can become a saddle point with $\sigma <1 $, or, otherwise, a saddle-focus).
A more delicate mechanism for destruction of pseudohyperbolicity is associated with the emergence of the so-called non-simple homoclinic tangencies, examples of which are shown in Fig.~\ref{htnonsimple}. In this case, as was established in \cite{Tat01,GGTat07,GOTat14}, stable periodic orbits, closed invariant curves and even nontrivial attracting invariant sets, e.g. small Lorenz-like attractors \cite{GOTat14}, can be born at bifurcations of such homoclinic tangencies.

From this fact we can draw an important conclusion for the theory of strange attractors of three-dimensional smooth maps: if such an attractor is genuine, then it should be either hyperbolic or pseudo-hyperbolic.
As for hyperbolic attractors, their mathematical theory is rather well developed, see e.g. \cite{An67,Sm67,AnSol91,GrinesandC,Kuzbook}.
We note, however, that, for a long time, the hyperbolic attractors had only purely mathematical interest, as no bifurcation mechanisms of their appearance in applications were known. The situation was changed after the papers by Turaev and Shilnikov \cite{TS95,TS97}, where it was shown that hyperbolic attractors (e.g. Smale-Williams and Anosov attractors) can be born at loss of stability of a stable periodic orbit under bifurcations of ``blue sky catastrophe'' type.
%\footnote{We note that such bifurcations can meet in models of neurodynamics \cite{GaASh}.}     Moreover, after a series of works by S. Kuznetsov
\cite{Kuz1,Kuz2,Kuz3,Kuz4}, it became known that such attractors are also found in physical models. Note that there are powerful analytical and computer methods for proving hyperbolicity of attractors.
% that to prove that the strange attractor in some model is hyperbolic, we have developed subtle qualitative and computer methods.
As we know, similar methods are now being created for detecting pseudohyperbolic attractors. %
In particular, in our recent works, new qualitative methods for investigating such attractors have been proposed, including such ones as new phenomenological bifurcation scenarios of appearance strange attractors of various types in one-parameter families \cite{GGS12,GGKT14,GG16,GK16}, and new search methods based on the effective using the so-called ``saddle charts'' \cite{GG16} and modified Lyapunov diagrams, new numerical methods of checking hyperbolicity etc. In the next sections we give some overview of these methods and give examples of some attractors which were found by these methods.
%obtained with the help of these methods,
We know that, at first examination, some of these attractors look to be genuine (pseudo-hyperbolic) ones.

\section{On phenomenological scenarios for the emergence of strange attractors of three-dimensional maps.}\label{sec:fen}

In this section we discuss questions of qualitative study of strange attractors of
three-dimensional maps. Moreover, the main attention will be paid to those attractors that can be pseudo-hyperbolic. Evidently, the spectrum $\Lambda_1,\Lambda_2,\Lambda_3$ of Luapunov exponents for orbits in such attractors must satisfy the necessary condition (\ref{lyapcond}), see also Remark~\ref{rem01}.
%Th for such attractors their
Besides,  we will restrict ourselves to the study of the so-called {\em homoclinic attractors}, i.e. such ones which contain only one fixed point of a map and its unstable manifold.

Here under attractor of a map $f$, following Ruelle \cite{Ruelle}, we mean {\em closed, invariant, stable, and chain-transitive set} ${\cal A}$.
As the stability, we will consider the usual asymptotic stability, which means that the attractor lies inside some absorbing domain ${\cal D}$, all points of which tend to
${\cal A}$ under positive iterations of the map $f$.
Recall that the chain transitivity, see e.g. \cite{AnBr85,TS98}, means that any two points on the attractor can be connected by an $\varepsilon$-orbit for any $\varepsilon>0$. The latter means that, for any two points $a, b \in {\cal A}$ and any $\varepsilon> 0$, in ${\cal A}$ there exist points $a = x_0, x_1, ..., x_ { N-1}, x_ {N} = b $, where $N = N(\varepsilon)$ is such that $x_i \in {\cal A}$ and $\mbox{dist} \; (x_{i + 1}, f (x_i)) <\varepsilon, \; i = 0, ..., N-1$.
The sequence of points $\{x_i\}$ is called an $\varepsilon$-orbit of the point $ x_0 $ of length $N + 1$, and the point $b$ is called $\varepsilon$-accessible from the point $a$.

Then we define the homoclinic attractor ${\cal A}$ with a fixed (periodic) point $O$ as a closed, invariant set consisting of points $\varepsilon$-accessible from the point $O$ for any $\varepsilon>0$. Hence, ${\cal A}$  is the prolongation of the point $O$.\footnote{That is, the closed invariant set containing all points $\varepsilon$-accessible from the point $O$. About the concept of prolongation in dynamical systems, see more details e.g. in \cite{AnBr85}.}

\begin{figure} [tb]
  \centerline{
  \includegraphics[width=12cm]{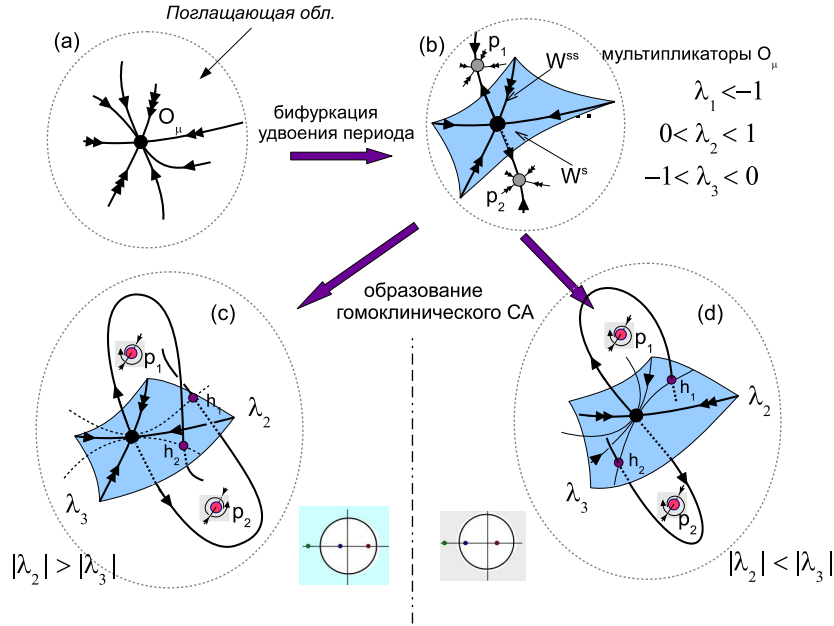}}
  \caption{{\footnotesize Schematic pictures of two phenomenological scenarios for the appearance of discrete homoclinic attractors: the case  of Lorenz-like attractor -- the path $(a)\to (b) \to (c)$;  the case of figure-eight attractor -- the path $(a) \to (b) \to (d) $.
}}
   \label{Lor8scen}
\end{figure}

In this case, geometrically the attractor ${\cal A}$, as a set in $R^3$, can be considered as the closure
(or, more exactly, the prolongation)
of the unstable manifold of its fixed point $O$.\footnote{This is true for the genuine attractor. However, we can not know this in advance. Nevertheless, such definition is agreed very well with computer study of attractors, when we can not see very small stable invariant subsets inside the attractor (e.g. stable periodic orbits of very large periods). Sometimes (e.g. for quasiattractors), these subsets can be visible and this correspond to the well-known phenomenon of appearance of ``windows of stability''. }
From this quite obvious observation, it can be concluded that geometrical as well as dynamical properties of the homoclinic attractor depend in much on its homoclinic structure, i.e. on a character of intersections of the stable and unstable invariant manifolds of the point $O$ itself. In this connection, in \cite{GGS12} we proposed rather simple phenomenological scenarios for the emergence of discrete homoclinic attractors of certain types  in one parameter families of maps starting with the simple attractor -- a stable fixed point. Two such scenarios are represented schematically in Figure~\ref{Lor8scen}.

We note two main features of these scenarios
%We note in these scenarios two main peculiarities, from our point of view.
The first is that, when the parameter changes, the stable fixed point $O$ loses stability under the supercritical period doubling bifurcation.
%The first is that the stable fixed point $O$ loses stability, when the parameter changes, under the supercritical period doubling bifurcation.
Immediately after this bifurcation, the point $O$ becomes a saddle with a one-dimensional unstable manifold, and in its neighbourhood a stable cycle $ (p_1, p_2) $ of period 2 is born (i.e. $f(p_1) = p_2$ and $f(p_2) = p_1 $), which is now the attractor.
Besides, the saddle point $O$ will have eigenvalues $\lambda_1,\lambda_2,\lambda_3$ such that $\lambda_1 <-1 $, $| \lambda_{2,3}| <1 $ and
$\lambda_2\lambda_3 <0$ (since $f$ is orientable). We assume that, at further change of parameter, the point $O$ no longer undergoes bifurcations and the cycle $ (p_1, p_2) $ loses its stability. How this happens does not yet matter, but what is important -- and this is the second main feature of these scenarios -- there is a global bifurcation associated with the creation of homoclinic intersections of the one-dimensional unstable $W^u$ and two-dimensional stable $W^s$ invariant manifolds of $O$. A configuration of these manifolds will be similar to what we see in Figure~\ref{Lor8scen} (c) and (d).

To explain how two such different configurations are formed, suppose, for definiteness, that $\lambda_2> 0$ and $\lambda_3 <0$ (here also
$\lambda_1 <-1$).  Then $W^u$ is partitioned by a point $O$ into two connected components, separatrices $W^{u +}$ and $W^{u -}$, which are invariant for $f^2$ and such that $f(W^{u +}) = W^{u -}$ and $f(W^{u-}) = W^{u +}$. For example, let $W^{u +}$ intersect $W^{s}_{loc}(O)$ at a point $h_1$, then $W^{u-}$ should intersect $W^{s}_{loc}(O)$ at the point $h_2 = f (h_1)$. The map $f$ in the restriction to $W^{s}_{loc} (O)$ is very simple: it has a stable fixed point $O$ of the type of a nonorientable node, since $\lambda_2 \lambda_3 <0 $.

In the case $|\lambda_2 |> |\lambda_3 |$, as in Fig.~\ref{Lor8scen}(c), in $W^{s}_{loc}(O)$ there exists a strongly stable invariant manifold
$W^{ss}$ that is an $f$-invariant curve,
tangent at the point $O$ to  the eigendirection corresponding to the negative multiplier $\lambda_3$.
The curve $W^{ss}$ splits the plane $W^{s}_{loc}(O)$ into two components $W^s_1$ and $W^s_2$. Since $\lambda_2> 0$ (and
$|\lambda_2 |> | \lambda_3 | $), each of these components is invariant under $f$, i.e. points from $W^s_1$ can not get into $W^s_2$ under iterations of $f$, and vice versa.
There is also a continuous family of smooth invariant curves on $W^{s}_{loc}$ that all enter the point $O$ touching the eigendirection corresponding to (positive) multiplier $\lambda_2$. Let the point $h_1$ belong to one of these curves, say $l_1$.
Then the curve $ l_2 = f (l_1) $ is also an invariant curve from this family, and $ h_2 \in l_2 $.
The curves $l_1$ and $l_2$ lie exactly in one component, either in $W^s_1$ or in $W^s_2$, and enter $O$, forming a ``zero-angle wedge'' configuration. Correspondingly, the configuration of unstable separatrices of the point $ O $, see Fig.~\ Ref {Lor8scen}(c), will resemble that typical for unstable separatrices of the Lorenz attractor. Therefore, the attractor arising here was named in \cite{GGOT13} ``Lorenz discrete attractor''.

Similar simple geometric arguments for the case $|\lambda_2| < |\lambda_3|$, as in Fig.~\ref{Lor8scen}(d), show that here the configuration of unstable separatrices of the point $O$ will be completely different. It is more like the configuration of separatrices in the attractor of the Poincar\'e map of a periodically perturbed two-dimensional system with a homoclinic figure-eight of saddle, \cite{GSV13}. Therefore, an attractor arising in this case was named in \cite{GGKT14} ``discrete figure-eight attractor'' (see Figs.~\ref{ExDLA1} and \ref{lorattr} which give an idea of the typical form of such an attractor).

We note that for both types of such attractors, the condition $\sigma> 1 $ (here $\sigma = |\lambda_1 \lambda_2 | $ in the ``Lorenz'' and
$\sigma = |\lambda_1 \lambda_3 | $ in the ``figure-eighth'' case, respectively) is very important, as it is  necessary in order to the attractor under consideration would be pseudohyperbolic.
Otherwise, these attractors are quasi-attractors of Lorenz or figure-eight type; or -- another possibility -- from a homoclinic configuration with $\sigma <1 $, a large enclosing it stable closed invariant curve (torus) can be created which further can be broken giving a rise of strange attractor a completely different nature (for example, ``torus-chaos''). Both these possibilities are well observed in computer experiments, see, for example, \cite{GGS12}.

These obvious observations tell us that, in the cases of saddle fixed points of other types, one can also expect the existence of homoclinic attractors whose configuration will depend essentially on eigenvalues of these points and, first of all, on their signs. In the case when there are complex conjugate eigenvalues, we can also expect the existence of discrete attractors of the spiral type, see \cite{GGS12,GGKT14} and part II of this paper.

\begin{remark}
{\rm However, our ``discrete'' Lorenz and figure-eight attractors differ substantially from their analogues obtained in Poincar\'e maps from periodically perturbed three-dimensional flows. Thus, for a small periodic perturbation of the system with the Lorenz attractor, we obtain a pseudo-hyperbolic attractor \cite{TS08} which has a saddle fixed point with all positive multipliers, and fixed points lie in the ``holes'' of the attractor.
In the case of the  Lorenz discrete attractor, the fixed point has two negative multipliers, and a period two orbit lies in ``holes''.  It seems that the discrete figure-eight attractors have no any flow analogues at all. This is due to the fact that if the corresponding system has a homoclinic figure-eight of saddle, then either this figure-eight is stable (an attractor) and then the attractor obtained will have
$\sigma <1 $, or it is unstable and then there is no attractor at all. This allows us to say that both the Lorenz discrete attractor  and the  discrete figure-eight attractor  are new.}
\label{rem1}
\end{remark}

The problem of {\em study and classification of homoclinic attractors for three-dimensional diffeomorphisms} itself was first stated in the paper \cite{GGS12}, although the first results on this subject were obtained in the paper \cite {GOST05}, in which discrete Lorenz attractors were found in three-dimensional H\'enon maps.
We note that a possibility of the appearance of such attractors at local bifurcations of triply degenerate fixed points, for example, having multipliers $ + 1; -1; -1 $, was investigated in the paper \cite{SST93}. Since the H\'enon map (\ref{3DHM1}) contains three parameters, such a point exists in it, and moreover, as shown in \cite{GOST05}, conditions from \cite{SST93} are fulfilled in this case. Thus, the main idea of our work \cite{GOST05} consisted in applying knowledge about the properties of degenerate local bifurcations to a specific situation. Obviously, this approach can also be used to the study of various models containing at least three parameters.

In \cite{GGS12} another idea was proposed based on the implementation of phenomenological scenarios of the appearance of strange homoclinic attractors possible in one-parameter families of three-dimensional maps.
Such scenarios as, for example, presented in Fig.~\ref{Lor8scen}, look quite realizable in specific systems and very convenient for computer research -- here, for example, you do not need to know all subtleties of global bifurcations leading to the appearance of homoclinic structures, but it is sufficient only to calculate/construct some basic simple characteristics (phase portrait, multipliers of the fixed point, Lyapunov exponents, etc). The idea of studying strange attractors by means of phenomenological scenarios involving two main bifurcation stages -- the loss of stability of a simple attractor (state of equilibrium, limit cycle, fixed point, etc.) and the appearance of a homoclinic attractor -- was first proposed by L.P. Shilnikov in the paper \cite {Sh86} where such a scenario was proposed to explain the phenomenon of emergence of spiral chaos in the case of multidimensional flows. We will discuss this scenario and its generalizations in the Part 2 of the paper.

Below we illustrate how these considerations can be used to study specific models
%Below
%In this paper
%we have illustrated how these ideas can be used for the study of strange homoclinic attractors in specific models.

\section{On numerical methods for study of pseudohyperbolic attractors.
} \label{sec:saddlemap}

The fact that the configuration of discrete homoclinic attractors
%under consideration
essentially depends on eigenvalues of their fixed points was used in the paper \cite{GG16} for the purposes of  classification of such attractors in the case of orientable three-dimensional maps.
If we restrict ourselves only to pseudohyperbolic homoclinic attractors, then such a classification problem turns out to be quite solvable if to distinguish
attractors by types of their homoclinic structures. In this case, as shown in \cite{GG16}, 5 different types of such pseudohyperbolic attractors are possible. All of them relate to the case when the fixed point is a saddle (all eigenvalues are real) with the one-dimensional unstable invariant manifold. Two of these types, the discrete Lorenz attractors and  discrete figure-eight attractors, can be observed in the case when the unstable eigenvalue
$\lambda_1$ is negative, i.e. $\lambda_1 <-1 $; and three other types of discrete attractors (the so-called ``double figure-eight'', ``super figure-eight'' and
``super Lorenz'' attractors) relate to the case when $\lambda_1> 1$, see \cite{GG16}.

To find such attractors in specific models, some fairly effective methods were proposed in \cite{GG16}.
%To find such attractors in specific models,
%In order to find such attractors in specific models,
One of them is the so-called ``method of saddle charts''.
%was proposed in \cite{GG16}.
We illustrate the essence of this method on the example of a {\em three-dimensional generalized H\'enon map} of the form
\begin{equation}
\bar x =y, \; \bar y = z,\; \bar z = Bx +Az+Cy  + f(y,z),
\label{GHM1-8}
\end{equation}
where
%$(x,y,z)\in \mathbb{R}^3$,
the nonlinearity $f$ depends only on coordinates $y$ and $z$ and, besides, $f(0,0)=0, f'_y(0,0)=f'_z(0,0)=0$.
The map (\ref{GHM1-8}) depends on three parameters $A, B$, and $C$ and has constant Jacobian equal to $B$. We will assume that $0<B<1$, i.e. the map is orientable and volume contracting.  Obviously,  any map of form
$\bar x = y, \; \bar y = z, \; \bar z = Bx + g (y,z)$ having a fixed point (for example, map (\ref{3DHM1}) for
$(1 + B-M_2)^2 + 4M_1> 0$) can be written in the form (\ref{GHM1-8}), if to move this point into the origin.

\begin{figure}[ht]
\centerline{\epsfig{file=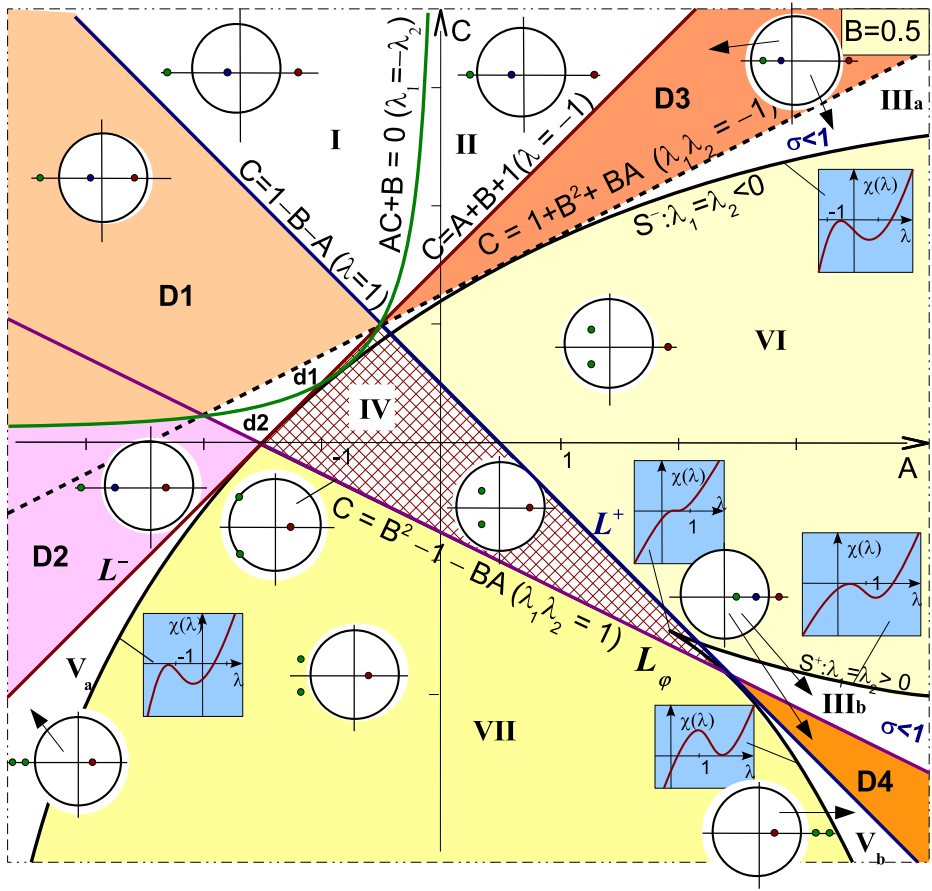, width=14cm
%, height=10cm
}}
\vspace{-1cm}
\caption{
{\footnotesize Examples of the saddle chart for map (\ref{GHM1-8}) with $B=0.5$.
}
}
\label{fig:charts05}
\end{figure}

The point $O(0,0,0)$ is a fixed point of map (\ref{GHM1-8}), the characteristic equation of (\ref{GHM1-8}) at this point has the form
\begin{equation}
\chi(\lambda)\equiv\lambda^3 - A \lambda^2 - C \lambda - B = 0.
\label{GHM3-8}
\end{equation}
It is important here that the eigenvalues of the point $O$ are functions of only the parameters $ A, B $ and $C$ and do not  depend on the nonlinearities $f(y,z)$. Then we can split the space of parameters $A$, $B$ and $c$ into domains corresponding to various types of location of eigenvalues of the point $O$ with respect to the unit circle.  We also distinguish the domains corresponding to  $\sigma> 1 $ and $\sigma <1 $ in the cases when the unstable manifold of $O$ are one-dimensional. Such a partition of the $(A,C)$-parameter plane for fixed $B$ ( a fixed value of the Jacobian) is called the {\em saddle chart} \cite{GG16}.\footnote{In the case of three-dimensional flows, a similar ``saddle chart'' for equilibrium states was proposed in \cite{book} in form  of a table, see there the  Appendix C.2.}
%We note that we also separate on this chart the domains corresponding to $\sigma> 1 $ and $\sigma <1 $ in the %case when the unstable manifold of the point $O$ is one-dimensional.
An example of such a saddle chart, with $B = 0.5$,
%and $ B = -0.5 $ are
is shown in Figure~\ref{fig:charts05}.

The domain \textsf{IV}, the so-called ``stability triangle''
(the domain $ \{C> B^2 -1 -BA \} \cap \{C <A + B + 1 \} \cap \{C <1-BA \} $) corresponds to the case when  the fixed point $O$ is asymptotically stable. For all other values of  $A$ and $C$ (except for bifurcation curves), the point $O$ is of saddle type -- it has eigenvalues both inside and outside the unit circle. Depending of the location of eigenvalues we select in the chart several domains shown in Fig.~\ref{fig:charts05}). The boundaries of the domains  consist of the  7 main curves. First, there are three bifurcation curves, with the values of the parameters for which the point $O$ has eigenvalues on the unit circle:
\begin{itemize}
\item
the curve $L_+$:  $\;\;C = 1 - B - A$ (when $\lambda = +1$);
\item
the curve $L_-$: $\;\;C = 1 + B + A$ (when $\lambda = -1$);
\item
the curve $L_\varphi$: $\;\;C = B^2 - 1 - BA$ at $-2 < A-B < 2$ (when $\lambda_{1,2}=e^{\pm i\varphi}$).
\end{itemize}
We note that the curve $ C = B ^ 2 - 1 - BA $ entirely belongs to the boundaries of the domains, but for $|A-B|>2$ it is not a bifurcation curve:
here  the point $O$  has eigenvalues  $(B, -|\lambda|, -|\lambda|^{-1})$  at  $ A-B  < -2 $  and  $(B, |\lambda|, |\lambda|^{-1})$ at  $ A - B > 2$.
Besides, the saddle chart contains 4 additional curves:
\begin{itemize}
\item
``the resonant curve'' $AC + B = 0, \;\; A <0$ (when $\lambda_1  = -\lambda_2$);
\end{itemize}
\begin{itemize}
\item
the curve ``$\sigma = 1$'' $C = 1 + B^2 + AB$ (when  $\lambda_1\lambda_2=-1$);
\end{itemize}
and
two curves of ``double roots''
\begin{itemize}
\item
$S^-$ (when  $\lambda_1  = \lambda_2 <0$),
\item
$S^+$ (when $\lambda_1  = \lambda_2 >0$).
\end{itemize}
The later two curves separate domains where $O$ has type of ``node'' and ``focus'' as well as  ``saddle'' and  ``saddle-focus''.
These curves have the following equations
$$
S^\pm: \;\;\left(\lambda_\pm\right)^3 - A \left(\lambda_\pm\right)^2 -C \lambda_\pm -B =0,
$$
where
$$
\displaystyle \lambda_\pm = \frac{A \pm \sqrt{A^2 + 3C}}{3}
$$
at $A^2 + 3C>0$ (i.e. $\lambda_\pm$ are the roots of equation $3\lambda^2 - 2A\lambda - C=0$).

In the saddle chart we remark especially the four regions D1, D2, D3 and D4, see  Fig.~\ref{fig:charts05}, where point $O$ is the saddle with one-dimensional unstable manifold, i.e., it has eigenvalues real multipliers $\lambda_1,\lambda_2,\lambda_3$ such that $|\lambda_1|>1, |\lambda_{2,3}| <1$, and such that the saddle value $\sigma = |\lambda_1|\cdot \max\{|\lambda_2|,|\lambda_3|\}$ is greater than 1. As we propose, only for values $A$ and $C$ from these domains, homoclinic attractors under consideration (containing the point $O$) can be pseudohyperbolic.
see  Fig.~\ref{fig:charts05}.
In other domains, except for the domain IV (the stability triangle), the point $O$ is either a saddle-focus, or a saddle with a two-dimensional unstable manifold, or has $\sigma <1 $. As we suppose, if the map has a homoclinic attractor (containing the point $O$) for the values of parameters from one of these last regions, then it is, in our opinion, a quasiatractor, see Remark~\ref{rem01}.

The use of saddle charts in numerical experiments  is a very convenient auxiliary tool in combining  with the Lyapunov diagram method. However, here we also modify the later method somewhat. Standardly, the Lyapunov diagram is a chart (in the $(A,C)$-parameter plane) consisting of colored areas corresponding to domains of parameters %eein it consists in constructing
%the spectrum  of Lyapunov exponents, in which
%the domains (in the $(A,C)$-parameter plane) painted in different colors correspond to regions of parameters
with different spectra of the  Lyapunov exponents $\Lambda_1> \Lambda_2>\Lambda_3$.
We use, in particular, the green color (in black and white drawings is also denoted by the number ``1'') -- for stable periodic regimes ($\Lambda_1 <0$); light blue color (number ``2'') -- for quasi-periodic regimes ($\Lambda_1 = 0$); yellow color (number ``3'') when $\Lambda_1> 0, \Lambda_2 <0 $, red color (number ``4'') when $\Lambda_1> 0, \Lambda_2 \sim 0 $, and blue color (number ``5''), when $\Lambda_1> \Lambda_2> 0 $ -- for strange
attractors.\footnote{The red domains (where $\Lambda_1> 0, \Lambda_2 \sim 0$) were especially highlighted in \cite{GOST05} -- these are such domains where the value of $\Lambda_2$ either always fluctuate very close to zero, or differ from zero by an amount (of the order of $10^{-5} $ or $10^{-6}$), comparable to the accuracy of  calculation of the exponents. Surprisingly, such domains turned out to be quite large, and this phenomenon (apparently related to the fact that the mapping on the attractor turned out to be very close to the time-discretization of some flow, for example, with the Lorenz attractor) was discussed in \cite{GOST05}.}  To these five colors, we added one more -- dark gray (number ``6'') to indicate {\em regions with homoclinic attractors}, when the numerically obtained points on the attractor approach the point $O$ by a very small distance (not less than $10^{- 4}$ for our numerics).

As was noted in Introduction, condition (2) for the spectrum $\Lambda_1,\Lambda_2,\Lambda_3$ of numerically obtained Lyapunov exponents  should be considered as one of necessary conditions for pseudohyperbolicity of strange attractors in three-dimensional maps. Moreover, this condition is evidently necessary for three-dimensional flows: here $\Lambda_2=0$ and, hence, $\Lambda_1>0$ automatically implies $\Lambda_1+\Lambda_2>0$; the inequality $\Lambda_1+\Lambda_2+\Lambda_3 <0$ follows from volume contracting properties of a flow near attractor.
However, it is well known that not all chaotic attractors for three-dimensional flows, and the more so for three-dimensional maps, are pseudohyperbolic. Thus, we strongly need  some additional numerical methods that give more confidence that the attractor is genuine. Such methods exist, for example, the Tucker methods of rigorous numerics \cite{Tucker99} based on interval arithmetic. However, the Tucker method is very difficult for us and too time-consuming for using it in our simple and standard numerics directed more for searching attractors then for their delicate studying. Instead this, we use a sufficiently simple but quite effective ``light method'' for verifying pseudohyperbolicity (LMP-method) of strange attractors of three-dimensional maps and flows (it can be applied even for four-dimensional flows, see Example 4 from Sec.5). This method was proposed in the paper \cite{KazTur17}.

The essence of the LMP-method consists in the fact that we take more attention for checking sufficient conditions for pseudohyperbolicity, see the Appendix.
In the case of three-dimensional maps (as well as flows), we consider condition (2) to be necessary one and assume that it holds. Then, the sufficient condition is related to the existence, at every point $x$ nearly attractor, of two transversal linear subspaces $N_1(x)$ and $N_2(x)$ such that
\begin{itemize}
\item[(i)] $\dim N_1 =1, \dim N_2 =2$;
\item[(ii)] $N_1(x)$ and $N_2(x)$ depend continuously on $x$;
\item[(iii)] $N_1(x)$ and $N_2(x)$ are invariant with respect to the differential $DT$ of the map $T$, i.e. $Df(N_1(x)) = N_1(f(x)),\;\;Df(N_2(x)) = N_2(f(x))$;
\item[(iv)] the map $T$ in the restriction to $N_1$ is uniformly contracting, and it in the restriction to $N_2$ extends exponentially two-dimensional volumes, and if in $N_2$ there is a contraction then it is uniformly weaker than the contraction in $N_1$.
\end{itemize}
In fact, these conditions are a simplified formulation of Definition~\ref{def:psevhyp} from the Appendix.

We note that the strongly contacting space $N_1(x)$ is one dimensional and it depends continuously on $x$. This means that angles $d\varphi$ between
any vectors $N_1(x)$ and $N_1(y)$ should be close for nearby $x$ and $y$. In fact, the LMP-method allows us to calculate these angles and, thus, to verify the continuity of the field $N_1$ of strong contacting directions at points of the attractor.
%Numerically,
%realization of
%the method is based on the well-known algorithm of calculation of the spectrum of Lyapunov exponents with some modifications.
The process of calculations consists of two stages. The first stage is standard: we calculate the spectrum of Lyapunov exponents $\Lambda_1,\Lambda_2,\Lambda_3$ (if condition (2) is not valid, we can stop calculations) and, in parallel, we store an array of data ${\cal N} = \{x_n\}$, where $x_{n+1} = f(x_n)$ and $n = 1,...k$, containing information about
%the coordinates of
points  $x_n$ on the attractor.
%and their sequences with respect to  iterations of the map.
The second step is not quite standard: we calculate the maximal Lyapunov exponent for backward iterations of the map using essentially the information obtained in the first stage. In particular, our backward iterations are forcibly attached to those points of the attractor that were obtained in the first stage.\footnote{Evidently, if we take any point on the attractor, then its backward iterations, sooner or later, depart far from the attractor. Thus, we can lose any information on the attractor. It is not the case when we take backward iterations exactly by the points on attractor taking these points from the first stage of the LMP-algorithm.} Note that the maximal Lyapunov exponent for backward iterations is equal with a minus sign to the minimal Lyapunov exponent $\Lambda_3$, and during these calculations we find vectors $N_1(x_n)$.
% for sufficiently large set $\{x_n\}, n=1,...,k,$ of points (where
%$k$ can be chosen so big as we want, depending on the problem).
As the final result of calculations, we construct the LMP-graph on the $(dx,d\varphi)$-coordinate plane, where $dx$ is the distance between two points $x$ and $y$ of the attractor and $d\varphi$ is the angle between vectors $N_1(x)$ and $N_1(y)$ (in fact, we construct the graph knowing points $x_i$ and $x_j$ and vectors $N_1(x_i)$ and $N_1(x_j)$ for all possible $i$ and $j$).

We note that if the attractor is pseudohyperbolic, which implies that the field $N_1(x)$ is continuous, the LMP-graph has to intersect the $d\varphi$-axis only at the origin
$(dx=0,d\varphi =0)$ or, if $N_1$ is not orientable, at the points $d\varphi =0$ and $d\varphi =\pi$. Thus, if the constructed  LMP-graph satisfies this property, we can conclude that our attractor should be surely pseudohyperbolic. On the other hand, if the  LMP-graph intersects  the $d\varphi$-axis in other points, except for $d\varphi =0$ and $d\varphi =\pi$ (or there is no a visible gap between the points of graph and the $d\varphi$-axis), we say that the attractor is a quasiattractor.

In section~\ref{Ex8Lor} we consider some examples three-dimensional maps and three- and four-dimensional flows with strange homoclinic attractors and verify these attractors for pseudohyperbolicity by the LMP-method.

\section{Examples of strange attractors in various models.}
  \label{Ex8Lor}

In this section we discuss several examples of strange attractors in three-dimensional maps, in a four-dimensional flow and in two nonholonomic models of rigid body dynamics. We consider only such attractors for which necessary conditions for pseudohyperbolicity (expressed by numerically obtained Lyapunov exponents) are satisfied. However, we show using the LMP-method that not all such attractors are genuine, some of them are, in fact, quasiattractors.

\subsection{Examples of homoclinic attractors in three-dimensional generalized H\'enon maps.}\label{3maps}

We note that there are many various methods to study chaotic dynamics in concrete models. One of the regular and reasonable approaches to this problem is related to the construction of diagrams of Lyapunov exponents.
%(or diagrams of dynamical regimes).
Namely in this way
%, using the modified diagram containing the ``dark grey spot'',
discrete Lorenz attractors were found in \cite{GOST05} for the three-dimensional H\'enon map of form (\ref{3DHM1}). Examples of such attractors are shown in Fig.~\ref{ExDLA1}. Now we can find such attractors, as they say, ``purposefully'', using our approach. To do this,
we consider the map (\ref{3DHM1}) in the following ``reduced to zero'' form
\begin{equation}
\bar x = y,\; \bar y = z,\; \bar z = B x + Az + Cy - z^2,
\label{3HM1priv}
\end{equation}
and take the saddle chart such as in Fig.~\ref {fig:charts05}(a) but constructed for the required fixed $B$, in our case for $B = 0.7$. Next, against the background of this chart, we numerically construct the modified 
\begin{figure}[ht]
\centerline{\epsfig{file=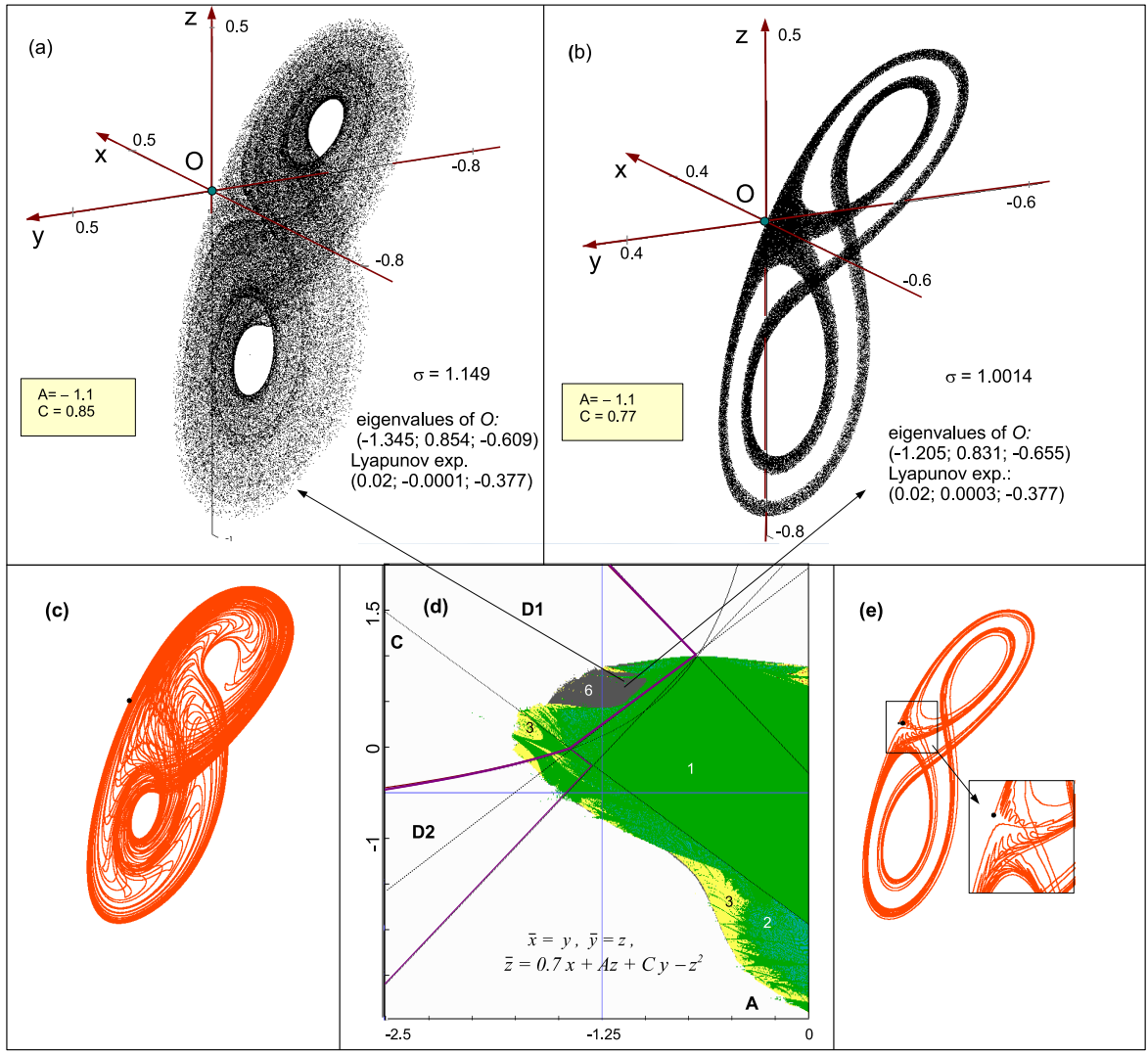, width=16cm
%, height=10cm
}}
\vspace{-0.5cm}
\caption{
{\footnotesize Discrete Lorenz attractor for map (\ref{3HM1priv}) with $B = 0.7$. }}
%In the middle: Diagram of Lyapunov exponents against the background of the ``saddle chart'' on the plane of parameters $A$ and $C$ for map (\ref{3HM1priv}) with $B = 0.7$. Right and left: phase portraits (containing approximately $10^4$ iterations of a single point) of attractors. }}
%(d) , (e)  numerically constructed one of the unstable separatrices of the point $O$ (other separatrix will look symmetrically, since $\lambda_1 <-1$), the presence of characteristic ``zigzags'' on the separatrix says that the attractor should be wild hyperbolic, i.e. it contains homoclinic tangencies; (f) magnification of  Fig. (e).}}
%
%
\label{lorattr}
\end{figure}
diagram of Lyapunov exponents. As a result, we get such a picture as in Fig.~\ref{lorattr}(d), where, in particular, the region of ``dark gray'' chaos intersects the region D1.  This suggests that, for the corresponding values of $A$ and $C$, a discrete Lorentz attractor can be  observed. The numeric results shown in Fig.~\ref{lorattr}(a) for $A=-1.1; C= 0.85$ and in Fig.~\ref{lorattr}(b) for  $A=-1.11; C= 0.77$, where we point out also values of the eigenvalues of $O$, the saddle value $\sigma$ of $O$, and the values of Lyapunov exponents $\Lambda_1,\Lambda_2,\Lambda_3$. In both cases we have that $\sigma>1$ and $\Lambda_1>0, \Lambda_1+\Lambda_2>0$. In Figs~\ref{lorattr}(c) and (e), we show also a behavior of one of the unstable separatrices of the point $O$ (the behavior of another separatrix is symmetric due to the unstable eigenvalue of $O$ is negative). We see that $W^u(O)$ has in both cases a homoclinic intersection with $W^s(O)$,
%and homoclinic tangencies are also
although, in the first case, typical zigzags in $W^u$ are seen more clearly than for the second case.

\begin{figure}[ht]
\centerline{\epsfig{file=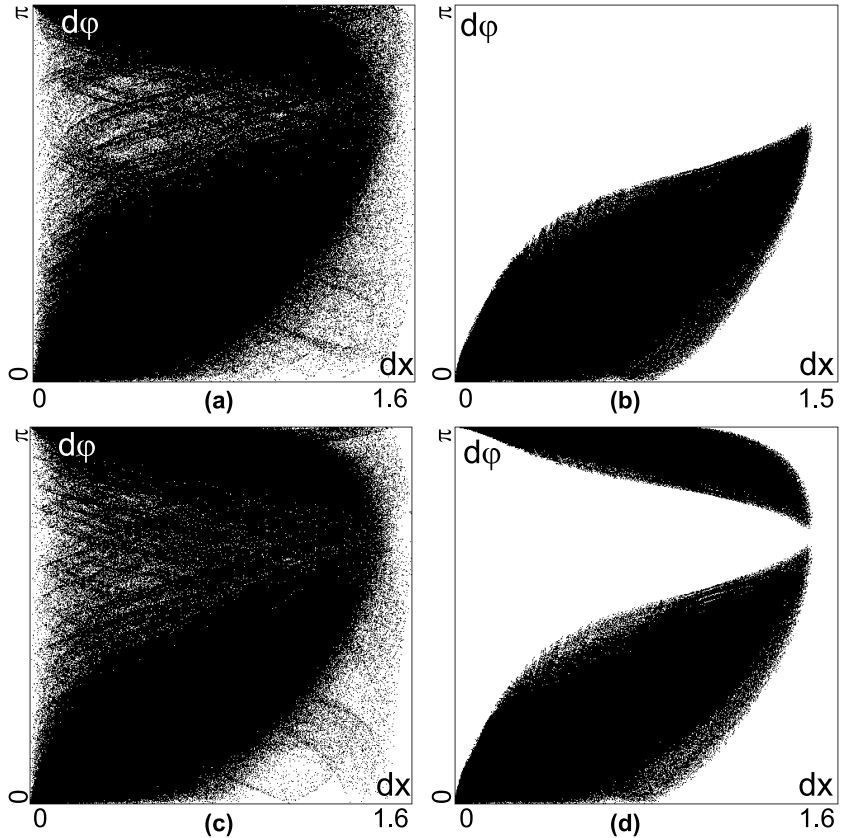, width=12cm}}
%, height=10cm
%}}
%\vspace{-1cm}
\caption{{\footnotesize LMP-graphs for attractors of Figs.~\ref{lorattr} (a) (left) and (b) (right).}}
\label{fig:DLorLMP}
\end{figure}

In Fig.~\ref{fig:DLorLMP} we represented the LMP-graphs for the discrete Lorenz attractors of Figs.~\ref{lorattr} (a) (left) and (b) (right).
The upper two figures show the LMP-graphs obtained by plotting every second iteration of the map, whereas,  in the lower figures the graphs are plotted for each iteration. In principle, there is no difference between figures (a) and (c), both of them look quite ``chaotic'' and show that the field $N_1$ in the case of attractors of Fig.~\ref{lorattr}(a) is not continuous. Thus, this attractor is certainly a quasiattractor.
It is not the case for figures (b) and (d). The evident difference between them can be explained by the fact that the the field $N_1$ of strong contracting directions is nonorientable here (this is inherited by the fact that a strongly stable eigenvalue of the fixed point $O$ is negative), and $N_1$ becomes orientable if to consider every second iteration of the map. Thus, we can conclude from the LMP-graphs of figures (b) and (d) that the attractor of Fig.~\ref{lorattr}(b) is pseudohyperbolic.
We can also remark that the LMP-graph of Fig.~\ref{fig:DLorLMP}(b), when every iteration is plotted, carries some important information about how close is our pseudo-hyperbolic 
attractor to its breaking down, when it becomes a quasiattractor. This can be estimated from the distance between the two ``whale'' components of the graph. When these components intersect the field $N_1$ immediately disappear and, as a consequence, nonsimple homoclinic tangencies like in Fig.~\ref{htnonsimple} appear.

\begin{figure}[ht]
\centerline{\epsfig{file=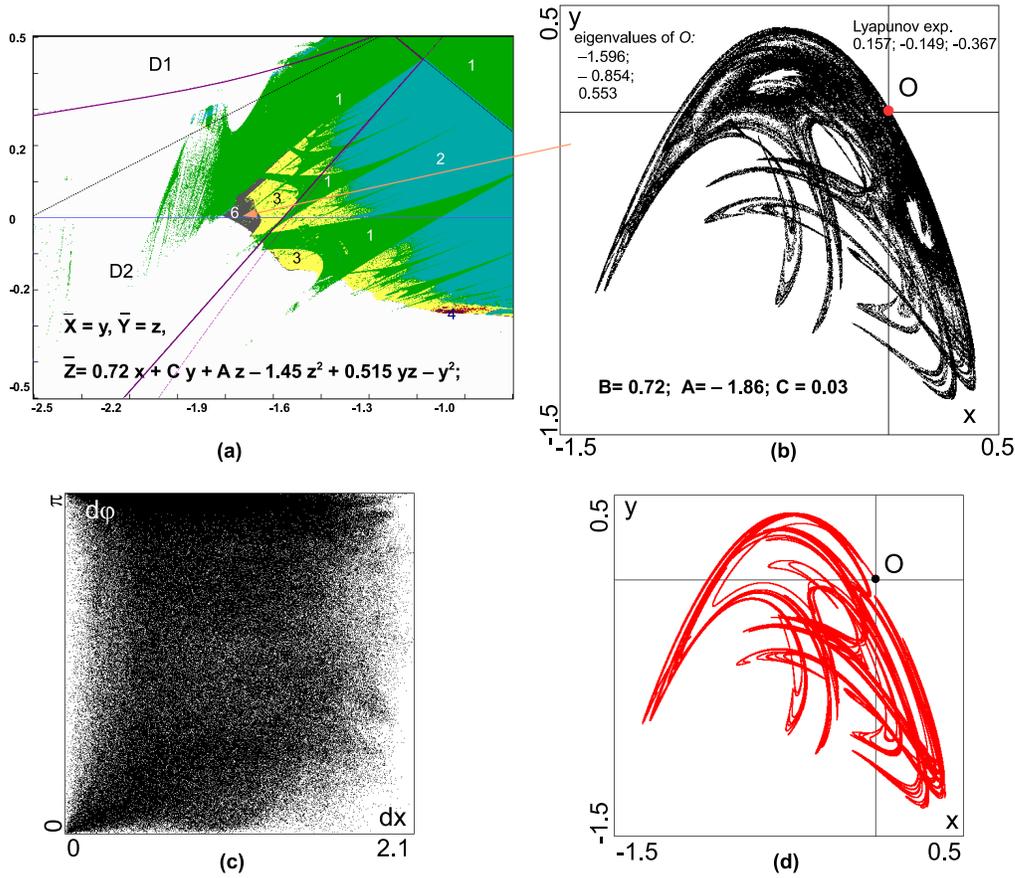, width=14cm
%, height=10cm
}}
%\vspace{-1cm}
\caption{
{\footnotesize An example of discrete figure-8 attractor in the corresponding generalized H\'enon map: (a) a fragment of the saddle chart on a background of the Lyapunov diagram;  (b) the projection of attractor on the $(x,y)$-plane; (c) the LMP-graph; (d) a behavior of one of the unstable separatrices of the fixed point $O$ is shown.
%Also, an additional information is given: the corresponding values of parameters $A$ and $C$, values of multipliers  and the saddle %value of $O$,  and the spectrum of Lyapunov exponents.
%
}
}
\label{Figure8}
\end{figure}

Obviously, in the case of the map (\ref{GHM1-8}), the saddle chart does not depend on the nonlinear terms $f(y,z)$. At the same time, the form of Lyapunov diagram on the $(A,C)$-parameter plane is  determined only  by these terms. With modern computers, the calculation of Lyapunov exponents
does not take much time, especially in the case of three-dimensional maps, and the  saddle chart for  map (\ref{GHM1-8} is constructed ``instantly''. This is all the more true when we are in the ``search stage''. In addition, as our experience shows, when varying nonlinearities one can see a certain tendency
\begin{figure}[ht]
\centerline{\epsfig{file=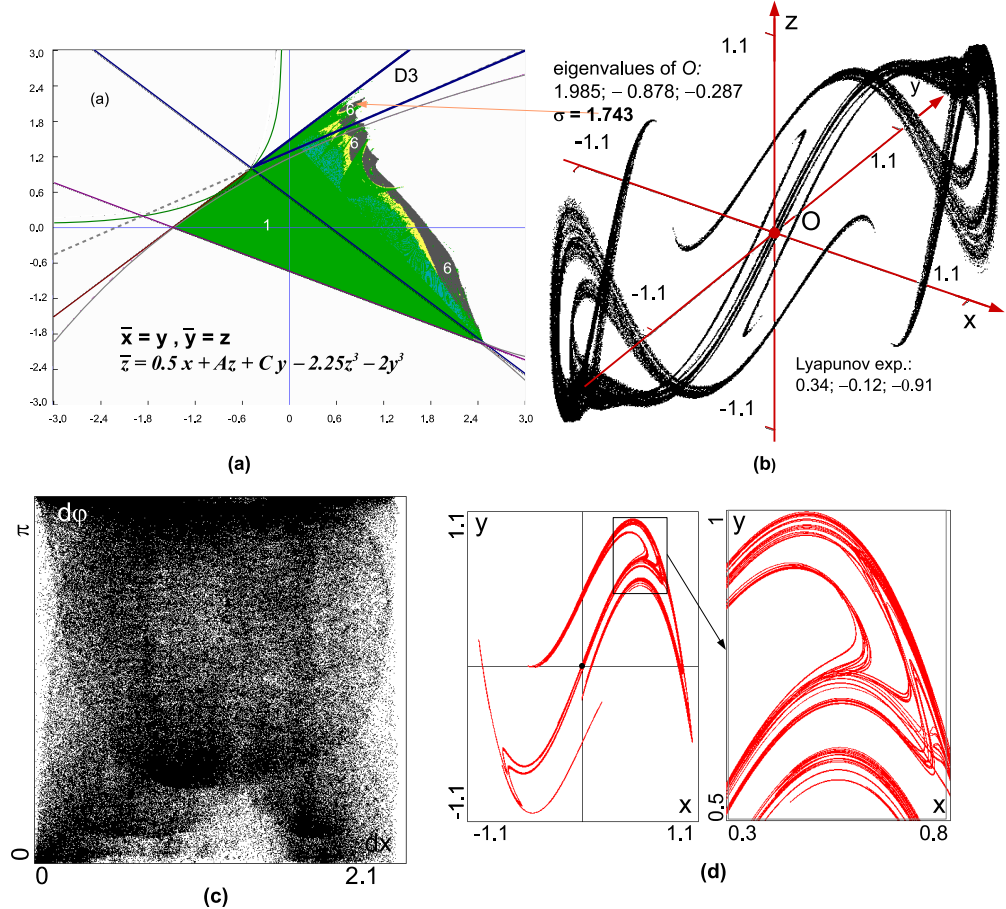, width=14cm
%, height=10cm
}}
\vspace{-1cm}
\caption{
{\footnotesize An example of discrete double figure-8 attractor in the corresponding generalized H\'enon map: (a) a fragment of the saddle chart on a background of the Lyapunov diagram;  (b) the phase portrait of attractor; (c) the LMP-graph; (d) a behavior of the unstable separatrices of the fixed point $O$ is shown (left) and a magnification of some fragment (right).
%Also, an additional information is given: the corresponding values of parameters $A$ and $C$, values of multipliers  and the saddle %value of $O$,  and the spectrum of Lyapunov exponents.
%
}
}
\label{DoubleFigure8}
\end{figure}
%
%, in particular,
in changing the location of the ``dark-gray spot'' (when the attractor is homoclinic). If desired, this spot can be  ``driven'' into any of the domains  of the saddle chart (except  for the ``stability triangle''), and, accordingly, we can find
the attractor of interest to us. By the same way, various homoclinic attractors of map (\ref{GHM1-8}) were found in \cite{GG16}. Some of them (for the values of  parameters $A$ and $C$ from the domains  D1, D2, D3 and D4) were presented as candidates for pseudohyperbolic attractors -- in particular, when the necessary condition (\ref{lyapcond}) is fulfilled for them.

In Fig.~\ref{Figure8} illustrations are shown that relate to the discrete figure-8 attractor of map (\ref{GHM1-8}) with the nonlinearity $f(y,z) = - 1.45 z^2 + 0.515 yz - y^2$ for values of parameters $B= 0.72;  A= - 1.86; C = 0.03$ belonging to the domain D2.  Although the necessary condition (\ref{lyapcond}) is fulfilled for the attractor, it looks as a typical quasiattractor that its LMP-graph of Fig.~\ref{Figure8}(c) confirms. Unfortunately, we can not find good examples of discrete figure-8 attractors in the case of three-dimensional H\'enon maps (however, we are sure that pseudohyperbolic attractors of such type exist here).

\begin{figure}[ht]
\centerline{\epsfig{file=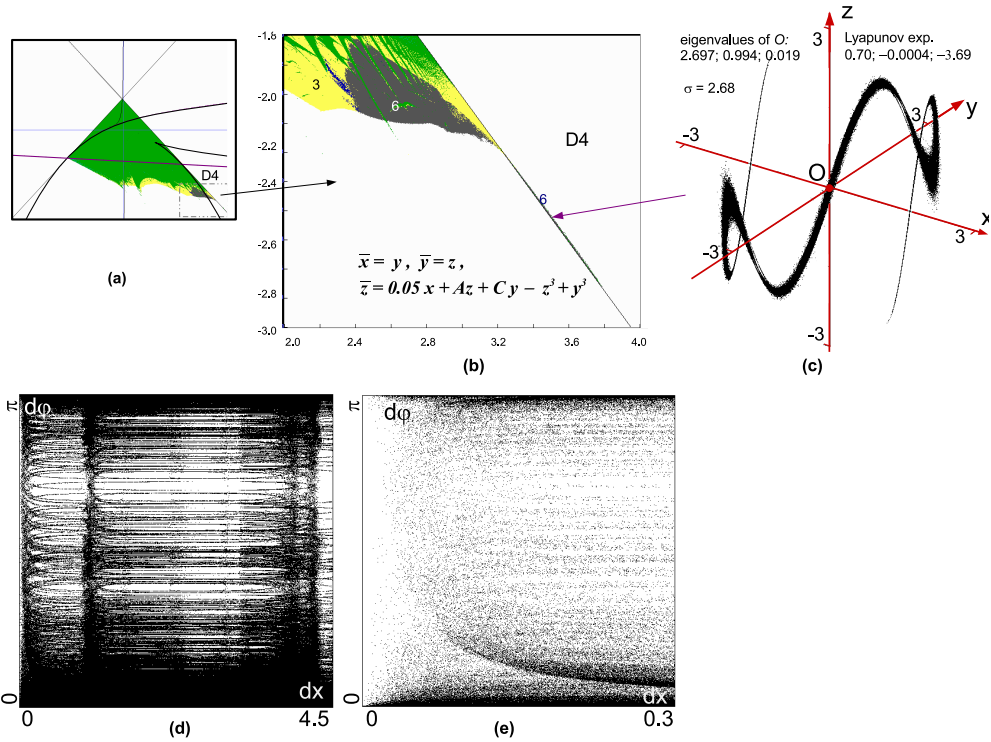, width=14cm
%, height=10cm
}}
%\vspace{-1cm}
\caption{
{\footnotesize An example of discrete super Lorenz attractor in the corresponding generalized H\'enon map: (a) a fragment of the saddle chart on a background of the Lyapunov diagram; (b) a magnification of some part of the Lyapunov diagram (a very thin dark-grey strip is seen in domain D4); (c) the phase portrait of attractor; (d) the LMP-graph; (e) a part of the LMP-graph near the axis $d\varphi$.
%Also, an additional information is given: the corresponding values of parameters $A$ and $C$, values of multipliers  and the saddle %value of $O$,  and the spectrum of Lyapunov exponents.
%
}
}
\label{SuperFigure8}
\end{figure}

In Fig.~\ref{DoubleFigure8} illustrations are shown that relate to the discrete double figure-8 attractor of map (\ref{GHM1-8}) with the cubic nonlinearity $f(y,z) = - 2.25z^3 - 2y^3$ for values of parameters $B= 0.5;  A=  0.82; C = 2.06$ belonging to the domain D3.  The necessary condition (\ref{lyapcond}) are satisfied again, but  LMP-graph of Fig.~\ref{DoubleFigure8}(c) shows that the attractor is certainly a quasiattractor.

In Fig.~\ref{SuperFigure8} illustrations related to one more homoclinic attractor, the so-called discrete super figure-8 attractor. This attractor is observed in map (\ref{GHM1-8}) with the cubic nonlinearity $f(y,z) = y^3 - z^3$ for values of parameters $B = 0.05; A = 3.702; C = -2.749$ belonging to the domain D4.  The necessary condition (\ref{lyapcond}) are satisfied again. Concerning the LMP-graph, see Fig.~\ref{SuperFigure8}(c), we have here very suspicious case, since nearly to the axis $dx=0$ we see a thin strip-like area where there are practically no points of the graph. This says more in favor of the fact that the attractor is pseudo-hyperbolic, although it exists in very small domain of parameters.

\subsection{Example of Shilnikov-Turaev wild spiral attractor.} \label{wildsfattr}

\begin{figure}%[ht]
\centerline{\epsfig{file=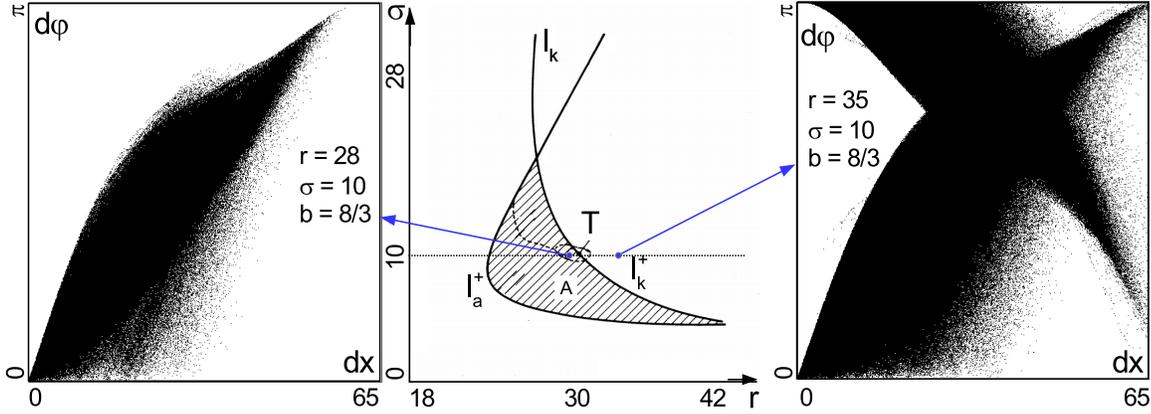, width=16cm
%, height=10cm
}}
\vspace{-0.5cm}
\caption{
{\footnotesize An illustration of results obtained by the LMP-method: (a) LMP-graph for the classical Lorenz values $\sigma = 10, r = 28, b = 8/3$; (b) the domain $A$ of the $(\sigma,r)$-parameter plane (for $b=8/3$) corresponding to the existence of the genuine Lorenz attractor -- this figure  is taken from \cite{BykovShilnikov92}; (c) LMP-graph for values $\sigma = 10, r = 35, b = 8/3$ from the right of the domain $A$. }}
\label{fig:LorenzPseudoTest}
\end{figure}

At first  we present, as a simple and illustrative example, the results obtained by the LMP-method for the classical Lorenz model
\begin{equation}
\left\{
\begin{array}{l}
\dot x = \sigma (y - x), \\
\dot y = x (r-z) - y, \\
\dot z = xy - bz.
\end{array}
\right.
\label{eq:Lorenz}
\end{equation}
It is well-known from the Tucker's paper \cite{Tucker99} that, for  the classical values of parameters
($\sigma = 10, r = 28, b = 8/3$), ``the Lorenz attractor exists'', i.e. it satisfies conditions of the Afraimovich-Bykov-Shilnikov geometrical model \cite{ABS77,ABS82}. In other words, it is pseudohyperbolic in terms of \cite{TS98}. Numerics from \cite{BykovShilnikov92} shows that this property holds for some region $A$, see Fig.~\ref{fig:LorenzPseudoTest}(b), of the $(\sigma,r)$-parameter plane (for $b=8/3$). The right boundary $l_k^+$ of $A$ corresponds to a violation of the conditions from \cite{ABS77,ABS82} and, as result, the attractor becomes a quasiattractor in the domain to the right of $l_k^+$.
In Fig.~\ref{fig:LorenzPseudoTest} there are also presented the LMP-graphs (a) for the classical values of the parameters  $\sigma = 10, r = 28, b = 8/3$ and for values  $\sigma = 10, r = 35, b = 8/3$ from the right of the line $l_k^+$. Thus, our LMP-test confirms that the first attractor is genuine, indeed. In the second case we see that the LMP-graph intersects
the axis $d\varphi$ at the points $d\varphi =0$ and $d\varphi =\pi$. However, the field $N_1(x)$ is orientable here as it should always be in the case of a flow. Thus, the second attractor  is certainly a quasiattractor.

\begin{figure}[ht]
\centerline{\epsfig{file=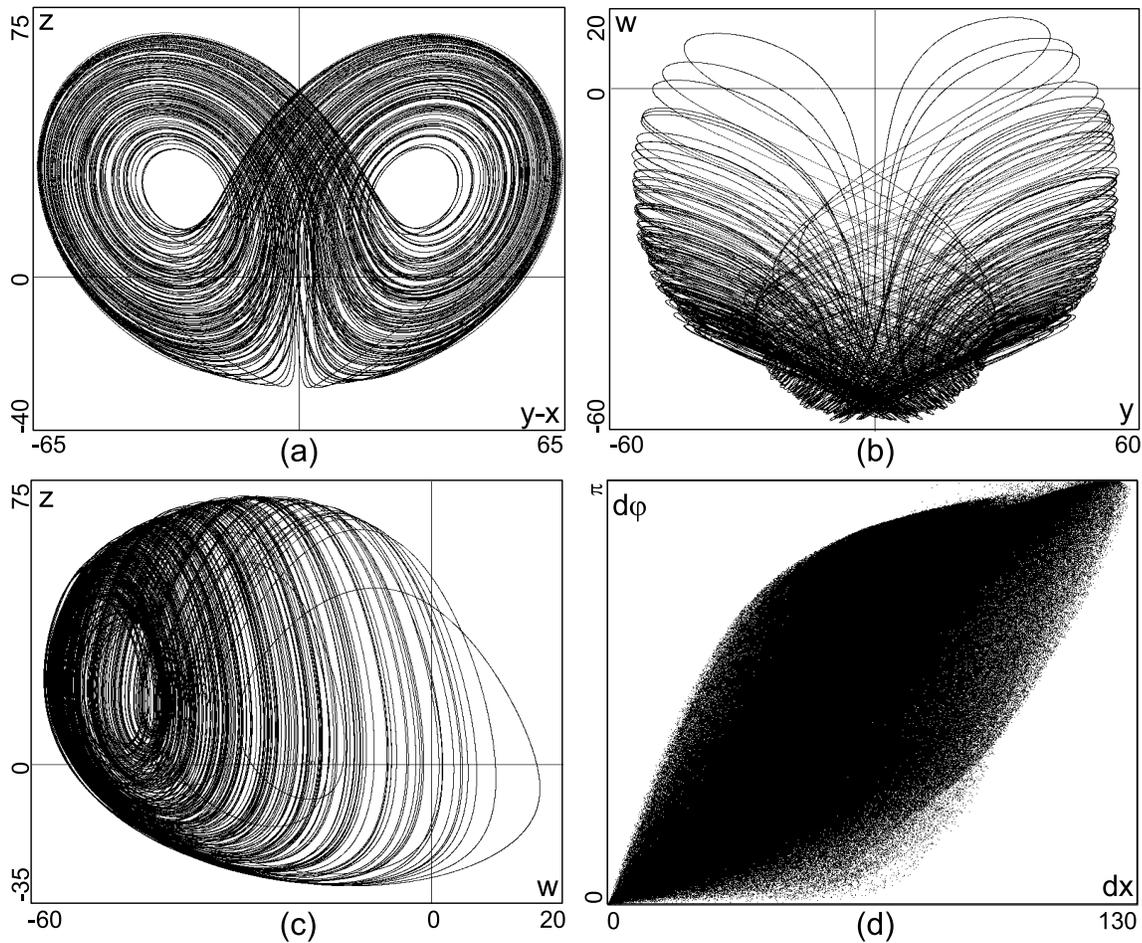, width=16cm
%, height=10cm
}}
\vspace{-0.5cm}
\caption{
{\footnotesize (a), (b) and (c) Projections of the attractor onto different two-dimensional planes;
(d) LMP-graph for the attractor. 
}
}
\label{wildattr}
\end{figure}

Now we consider another example of attractor found in \cite{KazTur17}, which is, in fact, wild spiral attractor, i.e. a pseudohyperbolic attractor containing an equilibrium of saddle-focus type.
Note that the base of theory of pseudohyperbolic attractors was laid in the paper \cite{TS98} by Shilnikov and Turaev, in which a geometric model of wild spiral attractor for a four-dimensional flow was constructed. This attractor contains a saddle-focus equilibrium $O$ with eigenvalues $\gamma, \lambda \pm i\omega, \tilde\lambda$, where $\gamma>0> \lambda > \tilde\lambda$ and, besides, $\gamma + 2\lambda >0$ and the divergence in $O$ is negative, i.e. $\gamma + 2\lambda + \tilde\lambda <0$. Thus, the point $O$ is pseudohyperbolic and
$\dim N_1(O) =1, \; \dim N_2(O) =3$; here the vector $N_1(O)$ is collinear to the eigenvector corresponding to the strong stable eigenvalue $\tilde\lambda$ and the three-dimensional plane $N_2(O)$ contains eigenvectors corresponding to three other eigenvectors of $O$ (thus, the plane $N_2(O)$ touches at $O$ the central unstable invariant manifold of the flow). The geometric model constructed in \cite{TS98} is in a sense similar to the Afraimovich-Bykov-Shilnikov model from \cite{ABS77,ABS82}, only the saddle is replaced by the saddle-focus and the flow under consideration is four-dimensional and, accordingly, the conditions of pseudohyperbolicity look to be more complicated. Until recently, the question on the existence of a concrete four-dimensional flow with wild hyperbolic attractor was open. In the paper \cite{KazTur17} examples of wild spiral attractors were found in some four-dimensional flows. One of such systems is the extended Lorenz system of the form
\begin{equation}
\left\{
\begin{array}{l}
\dot x = \sigma (y - x), \\
\dot y = x (r-z) - y, \\
\dot z = xy - bz + \mu w, \\
\dot w = -b w - \mu z,
\end{array}
\right.
\label{eq:LorenzModified}
\end{equation}
where $\sigma,r,b$ and $\mu$ are parameters. Note that at $\mu=0$ the system has an invariant three-dimensional plane $w=0$, in the restriction on which the system (\ref{eq:LorenzModified}) is the Lorenz system. When $\mu$ is nonzero this structure is broken and the Lorenz attractor existing e.g. at $\mu=0$ can evolute when $\mu$ varies. It would be quite interesting to track this evolution (e.g. when varying $\mu$ for fixed $\sigma,r,b$). However, we illustrate only one result from the paper \cite{KazTur17}, when  at
$$
r = 25, \sigma = 10, b = 8/3, \mu = 7
$$
the system (\ref{eq:LorenzModified}) has an attractor, whose projections onto two-dimensional planes (a) $\{w=0,x+y =0\}$; (b) $\{x=z=0\}$ and (c) $\{x=y=0\}$ are shown in Fig.\ref{wildattr}.
%and (d) its LMP-graph.
This attractor is spiral, since it contains the saddle focus equilibrium $O(0,0,0,0)$ with the eigenvalues
$$
\lambda_1 = \frac{1}{2}\left(\sqrt{(\sigma -1)^2 + 4\sigma r} -\sigma -1\right), \; \lambda_{2,3} = -b \pm i\mu, \; \lambda_4 = - \frac{1}{2}\left(\sqrt{(\sigma -1)^2 + 4\sigma r} +\sigma +1\right),
$$
i.e. $\lambda_1 = 10.93, \lambda_{2,3} = -8/3 \pm 7 i, \lambda_4 = -21.93$ for given values of parameters.
Thus, $O$ is a saddle-focus of type (3,1), i.e. with three-dimensional stable and one-dimensional unstable invariant manifolds, and, hence, $\dim N_1(O) =1, \; \dim N_2(O) =3$. %Also this attractor is surely pseudohyperbolic that one can see from its LMP-graph of Fig.\ref{wildattr}(d).
The necessary conditions
%Moreover, it was verified in \cite{KazTur2018} that the conditions
$\Lambda_1>0$, $\Lambda_1 + \Lambda_2 + \Lambda_3 > 0$ and $\Lambda_1 + \Lambda_2 + \Lambda_3 + \Lambda_4 < 0$
%, where $\Lambda_2 = 0$,
are also fulfilled here for numerically obtained Lyapunov exponents $\Lambda_1 = 2.19, \Lambda_2 = 0,  \Lambda_3 = -1.96, \Lambda_4 = -16.56$.
%for the first three Lyapunov exponents (and $L_4$ is strictly negative, so the sum of all four Lyapunov exponents is negative) and, besides,
Moreover, it is verified in \cite{KazTur17} that the field $N_1$ of strong stable directions at the points of attractor is continuous: the corresponding LMP-graph is shown in Fig.\ref{wildattr}(d). It is quite similar to that for the Lorenz attractor
%where the LMP-graph for the attractor is shown, the same as for the Lorenz attractor
(compare Figs.~\ref{wildattr}(d) and \ref{fig:LorenzPseudoTest}).

\subsection{Examples of pseudohyperbolic attractors in nonholonomic models of rigid body dynamics.} \label{nonmod}

In this section we review results on existence of strange attractors in two nonholonomic models of rigid body dynamics: the model of Celtic stone and the model of Chaplygin ball. Both models describe motions of a rigid body along the plane without slipping (this is a nonholonomic constrain) and with conserving the full energy. These nonholonomic and conservative approximations of the real motion of a rigid body along the plane allow to write the equations of its dynamics in the form of a five-dimensional system \cite{Markeev92} (three coordinates are momenta of impulse and two others are  two Euler angles) with the first integral that is the full energy $E$. Thus, in
the restriction on a level of the integral, the system becomes four-dimensional and we study its chaotic dynamics on the corresponding three-dimensional Poincar\`e section. Below we discuss the obtained results, for more details see \cite{GGK13,BKS14}.

Our first model is the Celtic stone model from \cite{GGK13}, see Section 4 there. This model depends on a lot of parameters characterizing physical and geometrical properties of the stone. As a natural control parameter we consider the value of the full energy $E$. When varying $E$ in an appropriate interval of its values we can observe a sequence of bifurcations in the corresponding one parameter family ${\cal F}_{E}$ of three-dimensional Poicar\'e maps leading from a stable fixed point to the discrete Lorenz attractor.

\begin{figure}[ht]
\centerline{\epsfig{file=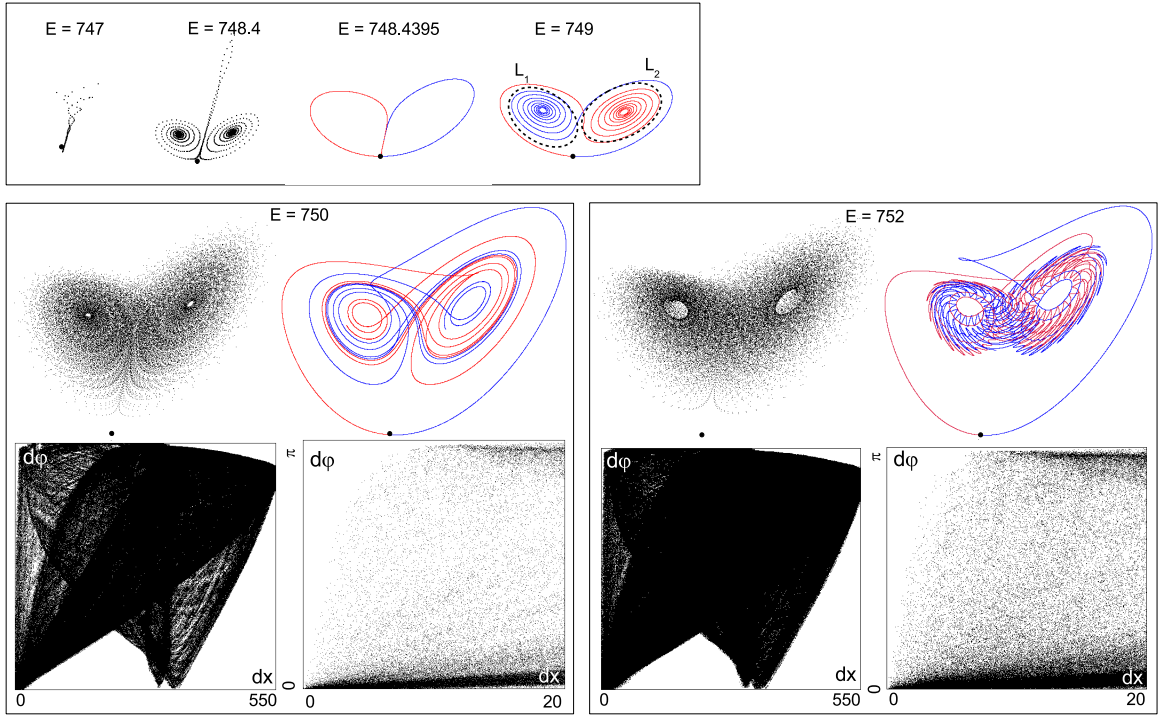, width=16cm
%, height=10cm
}}
\vspace{-0.5cm}
\caption{
{\footnotesize The upper line: first stages (from $E= 747$ to $E=749$) in developing dynamics to the discrete Lorenz attractors.
The middle line: $E=750$ and $E=752$, projections of attractors  (left) and unstable separatrices of $O$ (right) on some two-dimensional plane.
The bottom line: the corresponding LMP-graphs (left) and their magnifications (right).
}
}
\label{CeltLorenzFinal}
\end{figure}

The main stages of creation of such attractor are shown in Fig.~\ref{CeltLorenzFinal} when the parameter $E$ increasing from $E = 747$ to $E=E^* = 752$.
% $E = 747$ to $E=E^* = 750$.
%
At first, for $E < E_1 \simeq 747.61$,   the attractor is
a stable fixed point $O$, see Fig.~\ref{CeltLorenzFinal} for $E=747$. Then
it undergoes a period doubling bifurcation at $E = E_1$ and the stable cycle $P=(p_1,p_2)$
of period two becomes an attractor, see Fig.~\ref{CeltLorenzFinal} for $E=748.4$. Note that the point $O$ is now a saddle fixed point (with eigenvalues $\lambda_1<-1, 0<\lambda_2<1, -1<\lambda_3 <0$, where $|\lambda_3|<|\lambda_2|$ and $|\lambda_1||\lambda_2|>1$).
At $E=E_2 \simeq 748.4395$ a ``homoclinic figure-eight-butterfly'' of the unstable
manifolds (separatrices)
of the saddle $O$ is created, which gives rise then
to a saddle closed invariant curve $L = (L_1,L_2)$ of period two (where ${\cal F}_{E}(L_1)=L_2, {\cal F}_{E}(L_2)=L_1$),
the curve $L_1$ and $L_2$ surround the point $p_1$ and $p_2$, respectively, see Fig.~\ref{CeltLorenzFinal} for $E=749$.
At the same time, the unstable separatrices of $O$ are rebuilt and now, for
$E_2<E<E_3$, the left (the right) one wounds
to the right (the left) point of the cycle $P$. .
Moreover, together with the closed period-$2$ invariant curve $L$, the birth of an invariant limit
set $\Omega$ occurs here, \cite{Sh80}, which is not attracting yet.
As the numerical calculations show, for $E=E_3\sim 748.97$ the separatrices ``lie'' onto the stable
manifold of the curve $L$ and then leave it. Almost after that,
at $E = E_4 \sim 748.98$, the
period-$2$ cycle $P$ loses sharply the stability
under a subcritical torus birth bifurcation: the period-2 closed invariant curve $L$ merges with the cycle $P$
and after the cycle becomes of saddle-focus type.
% and the curve disappears.
The value of $E = E_4$ is the exact bifurcation moment of the discrete Lorenz attractor creation. We show in Fig.~\ref{CeltLorenzFinal} two examples of such attractors,  for $E=750$ and $E=752$.

Thus, this bifurcation scenario in the Celtic stone model fits into the overall scheme of appearance of discrete Lorenz attractors from Section~\ref{sec:fen}. However, it has a certain specific. So,
at the beginning, for $E$ close to $E_4$,  this attractor is quite unusual. Despite that it is a discrete attractor, unstable invariant manifolds of the point $O$ behave very similar to the flow case: here homoclinic intersections are invisible, since a size of splitting of the corresponding manifolds is comparable with the accuracy of calculations, see  Fig.~\ref{CeltLorenzFinal} for $E=750$. However, with increasing $E$ homoclinic intersection become more visible and typical zigzags appear in the unstable manifolds of $O$, see Fig.~\ref{CeltLorenzFinal} for $E=752$. Another specific is that the sequence of bifurcations when creating attractor is strikingly similar (one can say ``one-to-one'' for ${\cal F}_{E}^2$ ) to what happens in the Lorenz model \cite{Sh80}. Besides, this attractor is pseudohyperbolic that is seen from analysis of its LMP-graph, see Fig.~\ref{CeltLorenzFinal} where the LMP-graphs are shown for $E=750$ and $E=752$ both full-scale ones (left) and their enlarged fragments near the axis $d\varphi$. In the later we see that some neighborhoods of the axis $d\varphi$ do not contain any points of the graph. It follows surely that the attractors for $E=750$ and $E=752$ are pseudohyperbolic.

\begin{figure}[ht]
\centerline{\epsfig{file=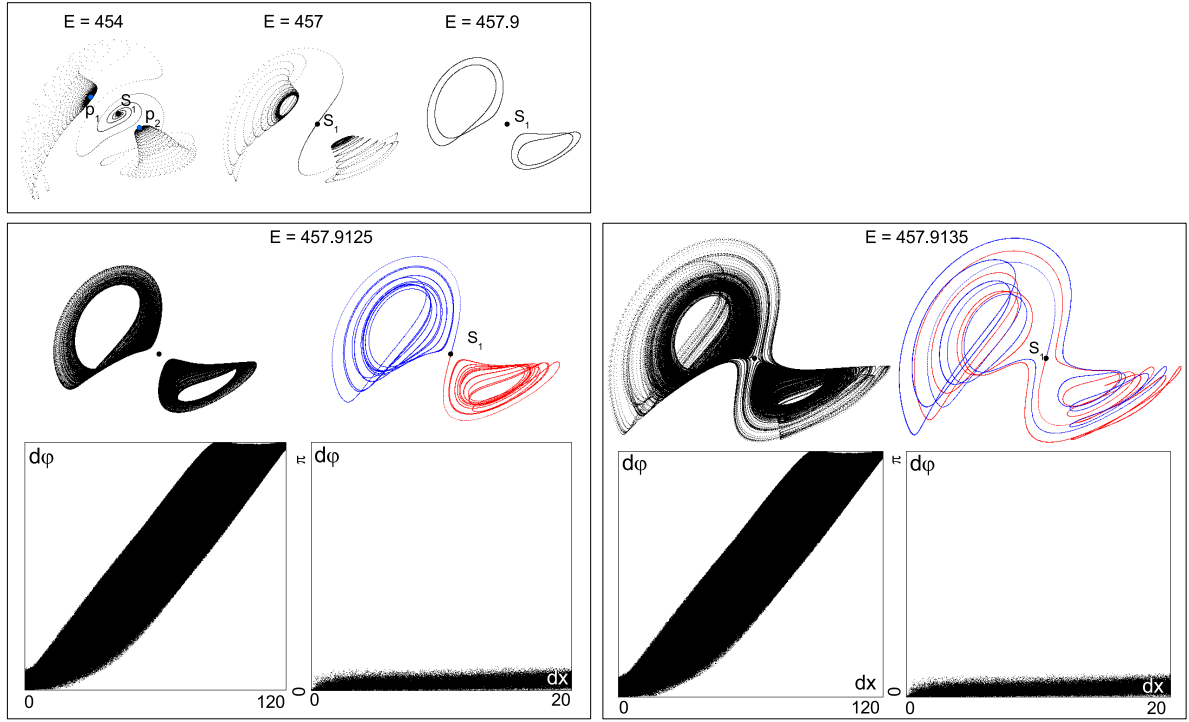, width=16cm
%, height=10cm
}}
\vspace{-0.5cm}
\caption{
{\footnotesize The upper line: first stages (from $E= 454$ to $E=457.9$) in developing chaotic dynamics.
The middle line: a two component torus-chaos for $E=457.9125$ and a discrete figure-eight attractor for $E=457.9135$; in both pictures there is shown  a two-dimensional projection for attractor (left) and for unstable separatrices of $S_1$ (right).
The bottom line: the corresponding LMP-graphs  and their magnifications.
}
}
\label{Figure8Final}
\end{figure}

The second our example is the nonholonomic model of Chaplygin ball from the paper \cite{BKS14}.
Figure~\ref{Figure8Final} shows the development of the attractor of the Poincar\'e map
in the model as the energy $E$ grows from $E= 454$.
At first, for $E_1 \simeq 417.5 < E < E_2 \simeq 455.95$, Fig.14(a),  the attractor is a period-2 orbit $(p_1,p_2)$ that emerges at $E=E_1$ along with a saddle period-2 orbit
$S=(s_1,s_2)$ as a result 
of a saddle-node bifurcation. Simultaneously, the system has a saddle fixed
point $S_1$: this point, a saddle-focus then a saddle, has a two-dimensional unstable manifold;
%then it has a two-dimensional unstable manifold;
then at
$E=E_3 \simeq 456.15$, the fixed point becomes a saddle with one-dimensional unstable manifold as a
result of a subcritical period-doubling bifurcation
when the saddle orbit $(s_1,s_2)$ merges with $S_1$. At
$E=E_2 \simeq 455.95$  the orbit $(p_1,p_2)$ loses the
stability at a supercritical Andronov-Hopf bifurcation and a stable period-2 closed invariant curve appears.
Thus, at $E > E_3$ the one-dimensional unstable separatrices of the saddle fixed point $S_1$ (with multipliers
$\lambda_1 < -1, |\lambda_{2,3}| <1$ and $\lambda_2\lambda_3 <0$ wind up
onto a stable closed invariant curve of period-2, see Fig.~\ref{Figure8Final} for $E=457$.
Next, several doublings of the invariant curve take
place, see Fig.~\ref{Figure8Final} for $E=459$, where a result of the first doubling is shown. The further growth of
$E$ leads to a strange attractor: at first a two component ``torus-chaos'', see Fig.~\ref{Figure8Final} for $E=457.9125$ and then a discrete figure-eight attractor, see Fig.~\ref{Figure8Final} for $E=457.9135$.

Note that at $E= 457.9135$, the fixed point $S_1$
has the eigenvalues $\lambda_1 \simeq -1.00907,\lambda_2 \simeq - 0.99732, \lambda_3 \simeq
0.98885$. Thus, the area-expansion conditions
$|\lambda_1\lambda_2| >1$ is fulfilled. Moreover, the Lyapunov
exponents for a random trajectory in the attractor are as follows: $\Lambda_1 \simeq 0.00063, \Lambda_2 \simeq - 0.00003, \Lambda_3 \simeq - 0.00492$, which gives $\Lambda_1+\Lambda_2>0$  and hints the pseudohyperbolicity. Besides, the corresponding LMP-graph in Fig.~\ref{Figure8Final} (bottom line) confirms this.
Note also that the ``torus-chaos'' of Fig.~\ref{Figure8Final} for $E=457.9125$ looks to be surely pseudohyperbolic that its
LMP-graph confirms.

\section{Appendix: to the definition of pseudo-hyperbolicity of maps.
} \label{App1}

We consider an $m$-dimensional diffeomorphism $f$. Let $Df$  be its differential.\footnote{Recall that the differential of the map
$f:\mathbb{R}^m \rightarrow \mathbb{R}^m$ in a point $x_0$  is a linear operator  $\displaystyle A = \frac{\partial f}{\partial x} \biggl|_{x = x_0} $ that maps a vector $\ell_{x_0}$ at $x_0$ to the vector $\ell_{x_1} = A\ell_{x_0}$ at the point $x_1 = f(x_0)$.}
An open domain ${\cal D} \subset\mathbb{R}^m$ is called {\em absorbing domain} of the diffeomorphism $f$ if
$f(\overline{{\cal D}}) \subset {\cal D}$.

\begin{df}
A diffeomorphism $f$ is called \textsf{pseudo-hyperbolic} on ${\cal D}$ if the following conditions are satisfied.
\begin{itemize}
\item[{\rm 1)}] Each point of ${\cal D}$ has two transversal linear subspaces $N_1$ and $N_2$,
which have complementary dimensions $(\dim N_1 = k \geq 1, \dim N_2 = m-k \geq 2)$, are continuously dependent  on the point and invariant under $Df$,
i.e
    $$
    Df(N_1(x)) = N_1(f(x)),\;\;Df(N_2(x)) = N_2(f(x)),
    $$
and such that, for every orbit $L:\{x_i \;|\; x_{i+1}= f(x_i),\; i=0,1,...; x_0\in {\cal D} \}$, its maximal Lyapunov exponent for
$Df\bigl|_{N_1}$
%corresponding to the subspace $N_1$
is strictly less than the minimal Lyapunov exponent for $Df\bigl|_{N_2}$, i.e. the following inequality holds:
\begin{equation}
\begin{array}{l}
\!\!\!\!\!\!\!\!\!\!\!\!\!\!\!\!\displaystyle \limsup\limits_{n\to\infty} \frac{1}{n}\;\ln\; (\!\!\sup\limits_{\begin{array}{c}u\in N_1(x_0)\\ \|u\|=1\end{array}}\!\! \|Df^n(x_0) u\|) < \\
%\qquad
\qquad\qquad
\displaystyle <\;\liminf\limits_{n\to\infty}\frac{1}{n}\;\ln \;(\!\!\inf\limits_{\begin{array}{c} v\in N_2(x_0)\\ \|v\|=1\end{array}}\!\! \|Df^n(x_0) v\|),
\end{array}
\label{eq:psh1}
\end{equation}
where $Df^n$ is an $(m\times m)$-matrix defined as
$$
Df^n = Df_{x_{n-1}} \cdot \ldots \cdot Df_{x_1} \cdot Df_{x_0},
$$
and $\limsup\limits_{n\to\infty}$ and $\liminf\limits_{n\to\infty}$ are the superior and inferior limits, respectively.

\item[{\rm 2)}]
 Diffeomorphism $f$ in the restriction to $N_1$ is uniformly contractive, that is, there exist constants $\lambda> 0$ and $C_1> 0 $ such that
\begin{equation}
\|D f^n(N_1)\| \leq  C_1 e^{-\lambda n}.
\label{eq:psh2}
\end{equation}

\item[{\rm 3})] Diffeomorphism $f$ in the restriction to $N_2$ extends exponentially $(m-k)$-dimensional volumes, that is, there exist constants
$\sigma> 0$ and $C_2>0$ such that\footnote{If $\dim N_2 = 1$, then the usual definition of uniform hyperbolicity is obtained, therefore we require in the definition that $\dim N_2 \ geq 2$.}
\begin{equation}
|\det D f^n(N_2) | \geq  C_2 e^{\sigma n}.
\label{eq:psh3}
\end{equation}
\end{itemize}
\label{def:psevhyp}
\end{df}

From Definition \ref{def:psevhyp} it immediately follows that:

\begin{itemize}
\item[$1^*$]
all orbits in ${\cal D}$ are unstable: each orbit has the positive maximal Lyapunov exponent
$$
\Lambda_{max}(x) = \limsup\limits_{n\to\infty} \frac{1}{n}\;\ln\;  \|Df^n(x)\| > 0.
$$
\end{itemize}

Note that the conditions of pseudohyperbolicity mean that whole ($(m-k)$-dimensional) volumes in $N_2$ are stretched under forward iterations of $f$.
% on $(m-k)$-dimensional subspaces of $N_2$.
This does not prohibit the existence of contraction directions in $N_2$,
%This does not prohibit that there may be compressive directions on $N_2$,
but any contractions along them should be uniformly weaker than any contraction in $N_1$. Thus, the uniform hyperbolicity can be considered as very specific case of pseudohyperbolicity, when all directions in $N_2$ are uniformly expanding, i.e. when  the inequality
%are weaker than the uniform hyperbolicity conditions -- the latter mean that
$\|Df^{-n}(N_2)\| <C e^{-\sigma n}$ holds.
%, i.e. the uniform expansion should take place in all directions in $N_2$.
Nevertheless, the same as in the case of hyperbolic systems, \cite{An67,TS98}, the following result is standardly proved here.

\begin{itemize}
\item[$2^*$]
The pseudohyperbolicity conditions are preserved for all sufficiently small $C^r$-perturbations of the system. Moreover, the spaces $N_1$ and $N_2$ vary continuously.
\end{itemize}

It follows from the statement $1^*$ that if a pseudohyperbolic diffeomorphism $f$ has an attractor in ${\cal D}$, then this attractor is strange and does not contain stable periodic orbits, which, as follows from the condition $2^*$, do not appear also for small smooth perturbations.
In other words, pseudo-hyperbolic attractors are genuine attractors.

\subsection*{Acknowledgement}    This paper was supported by grant 14-41-00044 of the RSF. Numerical experiments in Section~\ref{Ex8Lor} were supported by the RSF grant 14-12-00811. AG and SG thank
the RFBR (grants No. 16-01-00364, No. 16-51-10005 KO-a) and the Russian Ministry of Science and Education (project 1.3287.2017, target part) for supporting scientific researches.

\end{document}